\documentclass[11pt]{article}
\usepackage{amssymb}
\usepackage{amsmath}
\usepackage{amstext}
\usepackage{fancyhdr}
\usepackage{amsthm}

 \usepackage[all]{xy}

\input amssym.def
\input amssym.tex
\def\slfrac#1#2{\hbox{\kern.1em %
 \raise.5ex\hbox{\the\scriptfont0 #1}\kern-.11em %
 /\kern-.15em\lower.25ex\hbox{\the\scriptfont0 #2}}}

\newcommand{\hsp}{\hspace*{\parindent}}

\newcommand{\eeq}{\end{equation}}
\newcommand{\beql}[1]{\begin{equation}\label{#1}}

\makeatletter
\def\@sect#1#2#3#4#5#6[#7]#8{\ifnum #2>\c@secnumdepth
     \def\@svsec{}\else
     \refstepcounter{#1}\edef\@svsec{\csname the#1\endcsname.\hskip .75em }\fi
     \@tempskipa #5\relax
      \ifdim \@tempskipa>\z@
        \begingroup #6\relax
          \@hangfrom{\hskip #3\relax\@svsec}{\interlinepenalty \@M #8\par}%
        \endgroup
       \csname #1mark\endcsname{#7}\addcontentsline
         {toc}{#1}{\ifnum #2>\c@secnumdepth \else
                      \protect\numberline{\csname the#1\endcsname}\fi
                    #7}\else
        \def\@svsechd{#6\hskip #3\@svsec #8\csname #1mark\endcsname
                      {#7}\addcontentsline
                           {toc}{#1}{\ifnum #2>\c@secnumdepth \else
                             \protect\numberline{\csname the#1\endcsname}\fi
                       #7}}\fi
     \@xsect{#5}}
\def\@begintheorem#1#2{\it \trivlist \item[\hskip \labelsep{\bf #1\ #2.}]}

\def\plain{plain}\ifx\fmtname\plain\csname fi\endcsname
     
     \input docstrip
     \preamble

     Do not distribute the stripped version of this file.
     The checksum in the header refers to the documented version.

     \endpreamble
     \generateFile{here.sty}{t}{\from{here.doc}{}}
     \endinput
\fi
\ifcat a\noexpand @\let\next\relax\else\def\next{%
    \documentstyle[here,doc]{article}\MakePercentIgnore}\fi\next
\ifx\@Hxfloat\@Hundef\else\expandafter\endinput\fi
\let\@Hxfloat\@xfloat
\def\@xfloat#1[{\@ifnextchar{H}{\@HHfloat{#1}[}{\@Hxfloat{#1}[}}
\def\@HHfloat#1[H]{%
\expandafter\let\csname end#1\endcsname\end@Hfloat
\vskip\intextsep\vbox\bgroup\def\@captype{#1}\parindent\z@
\ignorespaces}

\def\end@Hfloat{\egroup\vskip \intextsep}

\makeatother

\catcode`\@=11

\renewcommand{\section}{%
        \setcounter{equation}{0}%
        \@startsection {section}{1}{\z@}{-3.5ex plus -1ex minus
        -.2ex}{2.3ex plus .2ex}{\large\bf}%
        }

\catcode`@=12

\begin{document}
\setlength{\baselineskip}{1.0\baselineskip}
\begin{center}
{\Large {\bf Representations and cohomologies of differential 3-Lie
algebras with any weight
 }}
 \\
\vspace{1.5\baselineskip} {\large{\bf Qinxiu Sun}}
\vspace{.2\baselineskip}

Department of Mathematics, Zhejiang University of Science and Technology, Hangzhou, 310023\\
E-mail: qxsun@126.com\\

\vspace{1.5\baselineskip} {\large{\bf Shan Chen }}
\vspace{.2\baselineskip}

Department of Mathematics, Zhejiang University of Science and Technology, Hangzhou, China, 310023\\
E-mail: 1844785776@qq.com\\

\end{center}

\setlength{\baselineskip}{1.25\baselineskip}

\begin{abstract}
The purpose of the present paper is to study
 cohomologies of differential 3-Lie algebras with any weight. We introduce the representation of a
differential 3-Lie algebra. Moreover, we develop cohomology theory
of a differential 3-Lie algebra. We also depict the relationship
between the cohomologies of a differential 3-Lie algebra and its
associated differential Leibniz algebra with weight zero. Formal
deformations, abelian extensions and skeletal differential 3-Lie
2-algebras are characterized in terms of cohomology groups

\end{abstract}

{\bf MR Subject Classification 2010}: 17B60, 17A30, 81R12.

\footnote{Keywords: differential 3-Lie algebra, representation,
 cohomology, $\mathcal{O}$-operator, abelian extension, 3-Lie 2-algebra,crossed module  }

\section{Introduction}
\hsp The notion of n-Lie algebra was introduced by Filippov in 1985
\cite{11}. These Lie algebras are presented with various names such
as Filippov algebra, Nambu-Lie algebra. The notion of n-Lie algebra
has close relationships with many fields in mathematics and
mathematical physics. For their applications, readers are referred
to \cite{16, 22, 30} and their references. 3-Lie algebras are
special types of n-Lie algebras and have close connection with the
supersymmetry and gauge symmetry of multiple coincident M2-branes
are applied to the study of the Bagger-Lambert theory \cite{5}.
Moreover, the n-Jacobi identity can be regarded as a generalized
Plucker relation in the physics literature.

The study on differential associative algebras has a long history
\cite{27} and it appears in various mathematical branches from
differential Galois theory, differential algebraic geometry to
differential algebraic groups etc. The algebraic abstraction of
difference equations led to a difference algebra which has been
considered largely in parallel to the differential algebra
\cite{13}. Naturally, $\lambda$-differential operator (differential
operator with any weight $\lambda$) was investigated in \cite{23},
which is a generalization of both the differential operator and the
difference operator. An algebra with $\lambda$-differential operator
is called a differential algebra with weight $\lambda$\cite{4}.

 Cohomologies and deformations are vital roles both in mathematics and physics.
 Algebraic deformation theory was studied by Gerstenhaber for rings and algebras in
\cite{20, 21}. Furthermore, Nijenhuis and Richardson considered the
deformation theory of Lie algebras \cite{17, 18}. More generally,
Balavoine \cite{2} investigated deformation theory of quadratic
operads. The cohomology theory and deformation theory for n-Lie
algebras was developed in \cite{1, 9, 19, 26, 31}.

The importance of differential algebras makes it compelling to
develop their cohomology theory with any weight. Cohomological and
deformation theories of differential Lie algebras with weight zero
were studied in \cite{32}. These results have been extended to the
case of associative algebras, Leibniz algebras, Pre-Lie algebras,
Lie triple systems and $n$-Lie algebras \cite{7, 8, 28, 29}. In
\cite{24}, the cohomologies, extensions and deformations of
differential algebras with any weight were developed.

In this paper, we would like to study the representation and
cohomology and deformation theory of differential 3-Lie algebras
with any weight.

The paper is arranged as follows. In Section 2, we introduce the
notion of a representation of differential 3-Lie algebras with any
weight. The matched pair of differential 3-Lie algebras are also
studied. In Section 3, we are devoted to the cohomology theory of
differential 3-Lie algebras. The relationship between the cohomogy
of differential 3-Lie algebras and the associated differential
Leibniz algebras with weight zero is also characterized. In Section
4, we discuss the deformations of differential 3-Lie algebras. We
also investigate the notions of Nijenhuis operators and
$\mathcal{O}$-operators on differential 3-Lie algebra $(\mathfrak g,
d_{\mathfrak g })$. Furthermore, we prove that $(V,d_{V})$ becomes a
new differential 3-Lie algebra and $(\mathfrak g, d_{\mathfrak g })$
is a representation of the new differential 3-Lie algebra
$(V,d_{V})$. Naturally, we study cohomologies of the new
differential 3-Lie algebra $(V,\varphi_{V})$. In section 5, we study
abelian extensions of differential 3-Lie algebras. Finally, we focus
on study skeletal and strict differential 3-Lie 2-algebras and
crossed modules using cohomologies.

Throughout the paper, let $k$ be a field of characteristic 0. Unless
otherwise specified, all vector spaces are over $k$.

\section{Representations of differential 3-Lie algebras}
\hsp A 3-Lie algebra is a vector space $\mathfrak g$ together with a
skew-symmetric multilinear map $[ \ , \ , \ ]:
 \wedge^{3}\mathfrak g\longrightarrow \mathfrak g$
satisfying
$$[x_1,x_2,[y_1,y_2,y_3]]=[[x_1,x_2,y_1],y_2,y_3]+[y_1,[x_1,x_2,y_2],y_3]+[y_1,y_2,[x_1,x_2,y_3]].\eqno (2.1)$$
for all $x_i,y_i\in \mathfrak g$.

A differential 3-Lie algebra with weight $\lambda$ is a 3-Lie
algebra $(\mathfrak g, [ \cdot, \cdot, \cdot])$ with a linear map
$d_{\mathfrak g} : \mathfrak g\longrightarrow \mathfrak g$, such
that for all $x_i\in \mathfrak g~(i=1,2,3,4,5)$, the following
holds:
\begin{eqnarray*}&&d_{\mathfrak g}([x_1, x_2, x_3]) \\&= &[d_{\mathfrak g}x_1, x_2, x_3] + [x_1, d_{\mathfrak g}x_2, x_3] + [x_1,
x_2, d_{\mathfrak g}x_3]+\lambda [d_{\mathfrak g}x_1, d_{\mathfrak
g}x_2, x_3] +\lambda [ x_1,d_{\mathfrak g}x_2, d_{\mathfrak g}x_3]
\\&&+\lambda [d_{\mathfrak g}x_1,x_2, d_{\mathfrak g}x_3] +\lambda^{2}
[d_{\mathfrak g}x_1,d_{\mathfrak g}x_2, d_{\mathfrak
g}x_3].~~~~~~~~~~~~~~~~~~~~~~~~~~~~~~~~~~~~~~~~~~~~~~~~~(2.2)
\end{eqnarray*}

{\bf Example 2.1.} Let $({\mathfrak g}, [ \ , \ , \ ])$ be
3-dimensional 3-Lie algebra with basis ${\varepsilon_1,
\varepsilon_2, \varepsilon_3}$ and the nonzero multiplication is
given by
$$[ \varepsilon_1 , \varepsilon_2 , \varepsilon_3 ]=\varepsilon_1.$$
Define linear map $d_{\mathfrak g}:{\mathfrak g}\longrightarrow
{\mathfrak g}$ by the matrix $\begin {bmatrix}
 d_{11}&d_{12}&d_{13} \\
d_{21}&d_{22}&d_{23} \\
d_{31}&d_{32}&d_{33}
\end {bmatrix}$
with respect to the basis ${\varepsilon_1, \varepsilon_2,
\varepsilon_3}$.

Then $({\mathfrak g}, [ \ , \ , \ ],d_{\mathfrak g})$ is a
differential 3-Lie algebra of weight $\lambda$ if and only if
\begin{eqnarray*}&&d_{\mathfrak g}([\varepsilon_1, \varepsilon_2, \varepsilon_3]) \\&= &[d_{\mathfrak g}\varepsilon_1, \varepsilon_2, \varepsilon_3]
+ [\varepsilon_1, d_{\mathfrak g}\varepsilon_2, \varepsilon_3] +
[\varepsilon_1, \varepsilon_2, d_{\mathfrak g}\varepsilon_3]+\lambda
[d_{\mathfrak g}\varepsilon_1, d_{\mathfrak g}\varepsilon_2,
\varepsilon_3] +\lambda [ \varepsilon_1,d_{\mathfrak
g}\varepsilon_2, d_{\mathfrak g}\varepsilon_3]
\\&&+\lambda [d_{\mathfrak g}\varepsilon_1,\varepsilon_2, d_{\mathfrak g}\varepsilon_3] +\lambda^{2}
[d_{\mathfrak g}\varepsilon_1,d_{\mathfrak g}\varepsilon_2,
d_{\mathfrak g}\varepsilon_3].
\end{eqnarray*}
By direct calculation, $({\mathfrak g}, [ \ , \ , \ ],d_{\mathfrak
g})$ is a differential 3-Lie algebra of weight $\lambda$ if and only
if $$d_{21}=d_{31}=0$$ and
$$d_{22}+d_{33}+\lambda(d_{11}d_{22}+d_{11}d_{33}+d_{22}d_{33})+\lambda^{2}(d_{11}d_{12}d_{33}-d_{11}d_{23}d_{32}-d_{13}d_{22}d_{33})=0.$$

For example, $({\mathfrak g}, [ \ , \ , \ ],d_{\mathfrak g})$ is
differential 3-Lie algebra of weight $\lambda=1$ with $$d_{\mathfrak
g}=\begin {bmatrix}
 1&1&2 \\
0&1&1 \\
0&4&1
\end {bmatrix}.$$

 In the light of \cite{9}, if $\mathfrak g$ is a 3-Lie algebra, then $L=\wedge^{2}\mathfrak
g$ is a Leibniz algebra with Leibniz bracket $[ \ ,\ ]_{F}$ given by
$$[ X ,Y ]_{F}=y_1\wedge [x_1,x_2,y_2]+[x_1,x_2,y_1]\wedge y_2
$$ for all $X=x_1\wedge x_2$ and
$Y=y_1\wedge y_2$. Through direct computation, we have\\

{\bf Proposition 2.2.} Let $(\mathfrak{g},d_{g})$ be a differential
$3$-Lie algebra of weight $\lambda$. Then $(L=\wedge^{2}\mathfrak g,
[ \ ,\ ]_{F}, d_{L})$ is a differential Leibniz algebra of weight
$\lambda$, where $d_{L}(x\wedge y)=d_{g}(x)\wedge y+x\wedge
d_{g}(y)+\lambda d_{g}(x)\wedge d_{g}(y)$ for any $x\wedge y\in L$.
\\

A representation of 3-Lie algebra $\mathfrak g$ consists of a vector
space $V$ together with a bilinear map $\rho:\mathfrak g\wedge
\mathfrak g\longrightarrow gl(V)$ such that for all $x_i\in
\mathfrak g~(i=1,2,3,4,5)$ the following are satisfied
$$\rho(x_1,x_2)\rho(x_3,x_4)=\rho([x_1,x_2,x_3],x_4)+\rho(x_3,[x_1,x_2,x_4])+\rho(x_3,x_4)\rho(x_1,x_2),\eqno (2.3)$$
$$\rho(x_1,[x_2,x_3,x_4])=\rho(x_3,x_4)\rho(x_1,x_2)-\rho(x_2,x_4)\rho(x_1,x_3)+\rho(x_2,x_3)\rho(x_1,x_4).\eqno (2.4)$$

{\bf Definition 2.3.} Let $(\mathfrak{g}, d_{g})$ be a differential
$3$-Lie algebra of weight $\lambda$. A representation of
$(\mathfrak{g}, d_{g})$ is a triple $(V,\rho, d_{V})$, where
$d_{V}\in {End}(V)$ and $(V,\rho)$ is a representation of $3$-Lie
algebra $\mathfrak{g}$, such that for any $x_1,x_2\in \mathfrak{g}$,
the following holds:
\begin{eqnarray*}d_{V}\rho(x_1,x_2)&=&\rho(d_\mathfrak{g}x_1, x_2)+\rho(x_1,
d_\mathfrak{g}x_2)+\lambda \rho(d_\mathfrak{g}x_1,
d_\mathfrak{g}x_2)+\rho(x_1,x_2)d_{V} \\&&+\lambda
(\rho(d_\mathfrak{g}x_1, x_2)+\rho(x_1, d_\mathfrak{g}x_2)+\lambda
\rho(d_\mathfrak{g}x_1,
d_\mathfrak{g}x_2))d_{V}.~~~~~~~~~~~~~~~~~~~~ (2.5)\end{eqnarray*}

{\bf Example 2.4.} Let $(\mathfrak g,d_{\mathfrak g})$ be a
differential $3$-Lie algebra with any weight $\lambda$ and define
bilinear map $ad:\mathfrak g\wedge \mathfrak g\longrightarrow
{gl}(\mathfrak g)$ by ${ad}(x_1,x_2)(y)=[x_1,x_2,y]$. Then
$(\mathfrak g,{ad},d_{\mathfrak g})$ is a representation of
differential $3$-Lie algebra $(\mathfrak g,d_{\mathfrak g})$, which
is called the adjoint
representation.\\

 {\bf Proposition 2.5.} Let $(\mathfrak{g_1},d_{\mathfrak{g_1}})$
and $(\mathfrak{g_2},d_{\mathfrak{g_2}})$ be two differential
$3$-Lie algebras, such that $(\mathfrak{g_1}, \mathfrak{g_2} , \rho,
\varrho)$ is a matched pair of 3-Lie algebras with bilinear maps
$\rho:\wedge^{2} \mathfrak{g_1}\longrightarrow {gl}(\mathfrak{g_2})$
and $\varrho:\wedge^{2} \mathfrak{g_2}\longrightarrow
{gl}(\mathfrak{g_1})$, where the 3-Lie algebra structure on
$\mathfrak{g_1}\bowtie\mathfrak{g_2}$ is given by
\begin{eqnarray*}[x_1+a_1,x_2+a_2,x_3+a_3]&=&[x_1,x_2,x_3]+\varrho(a_1,a_2)x_3+\varrho(a_3,a_1)x_2+\varrho(a_2,a_3)x_1\\&&+[a_1,a_2,a_3]
+\rho(x_1,x_2)a_3+\rho(x_3,x_1)a_2+\rho(x_2,x_3)a_1,
(2.6)\end{eqnarray*} for any $x_i\in\mathfrak{g_1},
a_i\in\mathfrak{g_2}~(i=1,2,3)$. Moreover, if $(\mathfrak{g_2},
\rho, d_{\mathfrak{g_2}})$ is a representation of differential
$3$-Lie algebra $(\mathfrak{g_1},d_{\mathfrak{g_1}})$ and
$(\mathfrak{g_1},\varrho,d_{\mathfrak{g_1}})$ is a representation of
differential $3$-Lie algebra $(\mathfrak{g_2},d_{\mathfrak{g_2}})$.
Define
$$(d_{\mathfrak{g_1}}+d_{\mathfrak{g_2}})(x+a)=d_{\mathfrak{g_1}}(x)+d_{\mathfrak{g_2}}(a)$$
for any $x\in d_{\mathfrak{g_1}},a\in d_{\mathfrak{g_2}}$.

Then $(\mathfrak{g_1}\bowtie
\mathfrak{g_2},d_{\mathfrak{g_1}}+d_{\mathfrak{g_2}})$ is a
differential $3$-Lie algebra. We call $(\mathfrak{g_1},
\mathfrak{g_2} ,d_{\mathfrak{g_1}},d_{\mathfrak{g_2}}, \rho,
\varrho)$ a matched pair of differential $3$-Lie algebras.
\begin{proof}
Using (2.5) and (2.6), we can verify it.
\end{proof}

More information on representation and matched pair of 3-Lie
 algebras, one can refer to \cite{3, 15}.\\

{\bf Proposition 2.6.} Let $(V,\rho,d_{V})$ be a representation of
 the differential $3$-Lie algebra $(\mathfrak{g},d_{g})$ with weight
$\lambda$. Define bilinear map
$\hat{\rho}:\wedge^{2}\mathfrak{g}\longrightarrow {gl}(V)$ by
$$\hat{\rho}(x_1,x_2)=\rho(x_1,x_2)+\lambda( \rho(d_\mathfrak{g}x_1, x_2)+
\rho(x_1, d_\mathfrak{g}x_2)+\lambda \rho(d_\mathfrak{g}x_1,
d_\mathfrak{g}x_2)),\eqno (2.7)$$  for any $x_1, x_2\in \mathfrak{g}
$. Then $(V,\hat{\rho},d_{V})$ is still a representation of
$(\mathfrak{g},d_{g})$. Denote it by $V_{\lambda}$.

\begin{proof} It is routine to verify that (2.3) and (2.4) hold for $\hat{\rho}$. We only check that (2.5) holds. In fact, for any
$x_1,
x_2\in \mathfrak{g}$, according to (2.5), we obtain
\begin{eqnarray*}&&\hat{\rho}(d_\mathfrak{g}x_1, x_2)+
\hat{\rho}(x_1, d_\mathfrak{g}x_2)+\lambda
\hat{\rho}(d_\mathfrak{g}x_1,
d_\mathfrak{g}x_2)+\hat{\rho}(x_1,x_2)d_{V}+\lambda(\hat{\rho}(d_\mathfrak{g}x_1,
x_2)\\&&+ \hat{\rho}(x_1, d_\mathfrak{g}x_2)+\lambda
\hat{\rho}(d_\mathfrak{g}x_1,
d_\mathfrak{g}x_2))d_{V}-d_{V}\hat{\rho}(x_1,x_2)
\\&=&\rho(d_\mathfrak{g}x_1,x_2)+\lambda\rho(d_\mathfrak{g}^{2}x_1,x_2)+\lambda\rho(d_\mathfrak{g}x_1,d_\mathfrak{g}x_2)
+\lambda^2\rho(d_\mathfrak{g}^{2}x_1,d_\mathfrak{g}x_2)\\&&+
\rho(x_1,d_\mathfrak{g}x_2)+\lambda\rho(x_1,d_\mathfrak{g}^{2}x_2)+\lambda\rho(d_\mathfrak{g}x_1,d_\mathfrak{g}x_2)
+\lambda^2\rho(d_\mathfrak{g}x_1,d_\mathfrak{g}^{2}x_2)\\&&+\lambda(
\rho(d_\mathfrak{g}x_1,d_\mathfrak{g}x_2)+\lambda\rho(d_\mathfrak{g}^{2}x_1,d_\mathfrak{g}x_2)+\lambda\rho(d_\mathfrak{g}x_1,d_\mathfrak{g}^{2}x_2)
\\&&+\lambda^2\rho(d_\mathfrak{g}^{2}x_1,d_\mathfrak{g}^{2}x_2))
+(\rho(x_1,x_2)+\rho(d_\mathfrak{g}(x_1),x_2)+\rho(x_1,d_\mathfrak{g}x_2)+\lambda\rho(d_\mathfrak{g}x_1,d_\mathfrak{g}x_2)d_{V}
\\&&
+\lambda(\rho(d_\mathfrak{g}x_1,x_2)+\lambda\rho(d_\mathfrak{g}^{2}x_1,x_2)+\lambda\rho(d_\mathfrak{g}x_1,d_\mathfrak{g}x_2)
+\lambda^2\rho(d_\mathfrak{g}^{2}(x_1),d_\mathfrak{g}x_2)\\&&+
\rho(x_1,d_\mathfrak{g}x_2)+\lambda\rho(x_1,d_\mathfrak{g}^{2}x_2)+\lambda\rho(d_\mathfrak{g}(x_1),d_\mathfrak{g}x_2)
+\lambda^2\rho(d_\mathfrak{g}x_1,d_\mathfrak{g}^{2}x_2)\\&&+\lambda
\rho(d_\mathfrak{g}x_1,d_\mathfrak{g}x_2)+\lambda^2\rho(d_\mathfrak{g}^{2}x_1,d_\mathfrak{g}x_2)+\lambda^2\rho(d_\mathfrak{g}x_1,d_\mathfrak{g}^{2}x_2)
+\lambda^3\rho(d_\mathfrak{g}^{2}x_1,d_\mathfrak{g}^{2}x_2))d_{V}\\&&-d_{V}(\rho(x_1,x_2)
+\rho(d_\mathfrak{g}x_1,x_2)+\rho(x_1,d_\mathfrak{g}x_2)+\lambda\rho(d_\mathfrak{g}x_1,d_\mathfrak{g}x_2))\\&=&0.
\end{eqnarray*}
\end{proof}

\section{Cohomology theory of differential 3-Lie algebras}
\hsp In this section, we will study cohomologies of a differential
3-Lie algebra with coefficients in its representation. We also
investigate the relation between the cohmologies of differential
3-Lie algebras and associated differential Leibniz algebra.
\subsection{Cohomologies of differential 3-Lie algebras}

\hsp Let us begin to recall the cohomology of $3$-Lie algebras
\cite{6, 31, 34}.

Let $\mathfrak g$ be a $3$-Lie algebra and $L=\wedge^{2}\mathfrak g$
be the associated Leibniz algebra. Suppose that $(\rho,V)$ is a
representation of $\mathfrak g$, the space
$C_{3-{\hbox{Lie}}}^{p}(\mathfrak g,V)$ of $p$-cochains ($p\geq 1$)
is the set of multilinear maps of the form
$$f:\wedge^{p-1}L\wedge\mathfrak g\longrightarrow V,$$
and the coboundary operator
$\partial:C_{3-{\hbox{Lie}}}^{p}(\mathfrak g,V)\longrightarrow
C_{3-{\hbox{Lie}}}^{p+1}(\mathfrak g,V)$ is as follows:
\begin{eqnarray*}&&(\partial f)(X_1,\cdot\cdot\cdot,X_{p},z)\\&=&\sum_{1\leq i<
k\leq {p}}
(-1)^{i}f(X_1,\cdot\cdot\cdot,\hat{X_i},\cdot\cdot\cdot,X_{k-1},[X_i,X_k]_F,X_{k+1},\cdot\cdot\cdot,X_{p}
,z)
\\&&+\sum_{i=1}^{p}(-1)^{i}f(X_1,\cdot\cdot\cdot,\hat{X_i},\cdot\cdot\cdot,X_{p},[X_i,z])+
\sum_{i=1}^{p}(-1)^{i+1}\rho(X_i)f(X_1,\cdot\cdot\cdot,\hat{X_i},\cdot\cdot\cdot,X_{p},z)\\&&
+(-1)^{p+1}(\rho(y_{p},z)f(X_1,\cdot\cdot\cdot,X_{p-1},x_{p})+\rho(z,x_{p})f(X_1,\cdot\cdot\cdot,X_{p-1},y_{p}))~~~~~~~~~~~~~(3.1)
\end{eqnarray*}
for all $X_i=x_{i}\wedge y_{i}\in L $ and $z\in \mathfrak g$. Denote
the corresponding cohomology group by $H^{*}_{3-Lie}(\mathfrak g,
V)$.

Moreover, let $(\mathfrak{g},d_{g})$ be a differential $3$-Lie
algebra with weight $\lambda$ and $(V,\rho,d_{V})$ be its
representation. Then $V_{\lambda}$ is a new representation of
$(\mathfrak{g},d_{g})$, which is mentioned in Proposition 2.6.

Denote $C_{3-{\hbox{DLie}}}^{p}(d_\mathfrak
g,d_V)=\hbox{Hom}(\wedge^{p-1}L\wedge\mathfrak g,V)$ by the set of
$p$-cochains ($p\geq 1$) of the differential operator $d_{\mathfrak
g}$ with coefficients in representation $(V,\hat{\rho},d_{V})$, and
the coboundary operator
$$\partial_\lambda:C_{3-{\hbox{DLie}}}^{p}(d_\mathfrak
g,d_V)\longrightarrow C_{3-{\hbox{DLie}}}^{p+1}(d_\mathfrak g,d_V)$$
is given by
\begin{eqnarray*}&&(\partial_\lambda f)(X_1,\cdot\cdot\cdot,X_{p},z)\\&=&\sum_{1\leq i<
k\leq {p}}
(-1)^{i}f(X_1,\cdot\cdot\cdot,\hat{X_i},\cdot\cdot\cdot,X_{k-1},[X_i,X_k]_F,X_{k+1},\cdot\cdot\cdot,X_{p}
,z)
\\&&+\sum_{i=1}^{p}(-1)^{i}f(X_1,\cdot\cdot\cdot,\hat{X_i},\cdot\cdot\cdot,X_{p},[X_i,z])+
\sum_{i=1}^{p}(-1)^{i+1}\hat{\rho}(X_i)f(X_1,\cdot\cdot\cdot,\hat{X_i},\cdot\cdot\cdot,X_{p},z)\\&&
+(-1)^{p+1}(\hat{\rho}(y_{p},z)f(X_1,\cdot\cdot\cdot,X_{p-1},x_{p})+\hat{\rho}(z,x_{p})f(X_1,\cdot\cdot\cdot,X_{p-1},y_{p}))~~~~~~~~~~~~~(3.2)
\end{eqnarray*}
for any $f\in C_{3-{\hbox{DLie}}}^{p}(d_\mathfrak g,d_V)$ and
$X_1,\cdot\cdot\cdot,X_{p}\in L,~z\in \mathfrak g$.\\

{\bf Definition 3.1.} The cohomology of the cochain complex
$(C_{3-{\hbox{DLie}}}^{*}(d_\mathfrak g,d_V),\partial_\lambda)$,
denoted by $H_{DO_\lambda}^{*}(d_{\mathfrak g},d_{V})$, is called
the cohomology of the differential operator $d_{\mathfrak g}$ with
coefficients in the representation
$V_{\lambda}$.\\

{\bf Remark 3.2.} $C_{3-{\hbox{DLie}}}^{p}(d_\mathfrak
g,d_V)=C_{3-{\hbox{Lie}}}^{p}(\mathfrak g,V)$ as vector spaces, but
they are not
equal as cochain complexes except $\lambda=0$.\\

{\bf Proposition 3.3.} Define linear map
$\delta:C_{3-{\hbox{Lie}}}^{p}(\mathfrak g,V)\longrightarrow
C_{3-{\hbox{DLie}}}^{p}(d_\mathfrak g,d_V)$ by
$$\delta f(X_1,\cdot\cdot\cdot,X_{p},z)=
\sum_{k=1}^{2p+1}\lambda^{k-1}f^{k}(X_1,
\cdot\cdot\cdot,X_{p},z)-\varphi_{V}f(X_1, \cdot\cdot\cdot,X_{p},z),
$$
where $$f^{k} = f(\underbrace{I,\cdot\cdot\cdot , d_{\mathfrak
g},\cdot\cdot\cdot , I}_{d_{\mathfrak g}~~ \hbox{appears} ~~k
~~\hbox{times}})$$ and $I$ is the identity map on $\mathfrak g$.
Then $\delta$ is a cochain map from the cochain complex
$(C_{3-{\hbox{Lie}}}^{*}(\mathfrak g,V),\partial)$ to
$(C_{3-{\hbox{DLie}}}^{*}(d_\mathfrak g,d_V),\partial_\lambda)$,
that is, the following diagram is commutative:
$$\xymatrix{
 C_{3-{\hbox{Lie}}}^{n}(\mathfrak g,V) \ar[d]_{\partial} \ar[r]^{\delta} & C_{3-{\hbox{DLie}}}^{n}(d_\mathfrak g,d_V) \ar[d]^{\partial_{\lambda}} \\
  C_{3-{\hbox{Lie}}}^{n+1}(\mathfrak g,V)\ar[r]^{\delta} & C_{3-{\hbox{DLie}}}^{n+1}(d_\mathfrak g,d_V)
  .}$$

\begin{proof}
See Appendix.
\end{proof}

Denote the set of $p$-cochains by
$$C_{3-D_\lambda}^{p}(\mathfrak g,V)=:C_{3-{\hbox{Lie}}}^{p}(\mathfrak g,V)\times
C_{3-{\hbox{DLie}}}^{p-1}(d_\mathfrak g,d_V),~~p\geq 2$$ and
$$C_{3-D_\lambda}^{1}(\mathfrak g,V)=\hbox{Hom}(\mathfrak g,V).$$

Define the linear map
$\partial_{D_\lambda}:C_{3-D_\lambda}^{p}(\mathfrak
g,V)\longrightarrow C_{3-D_\lambda}^{p+1}(\mathfrak g,V)$ by
$$\partial_{D_\lambda}(f,g)=(\partial f,\partial_\lambda
g+(-1)^{p}\delta f)$$ for any $(f,g)\in
C_{3-D_\lambda}^{p}(\mathfrak g,V)$.

In view of Proposition 3.3, we have \\

{\bf Theorem 3.4.} The operator $\partial_{D_\lambda}$ is a
coboundary operator, that is,
$\partial_{D_\lambda}\partial_{D_\lambda}=0$.\\

Associated to the representation $(V,\rho,d_{V})$, we obtain a
cochain complex $(C_{3-D_\lambda}^{*}(\mathfrak g,V),
\partial_{D_\lambda})$. Denote the cohomology group of this cochain
complex by $H^{*}_{3-D_\lambda}(\mathfrak g, V)$, which is called
the cohomology of the differential $3$-Lie algebra
$(\mathfrak{g},d_{g})$ with coefficients in the representation
$(V,\rho,d_{V})$.

In view of the definitions of $C_{3-{\hbox{DLie}}}^{p}(d_\mathfrak
g,d_V)$ and $C_{3-D_\lambda}^{*}(\mathfrak g,V)$, we easily get an
short exact sequence of cochain complexes,
$$0\longrightarrow C_{3-{\hbox{DLie}}}^{n-1}(d_\mathfrak
g,d_V) \longrightarrow C_{3-D_\lambda}^{n}(\mathfrak g,V)
\longrightarrow C_{3-Lie}^{n}(\mathfrak g,V)\longrightarrow 0.$$

Moreover, combining the Snake Lemma, we have the following long
exact sequence of cohomology groups

$$\cdot\cdot\cdot\longrightarrow H_{DO_\lambda}^{n-1}(d_{\mathfrak
g},d_{V})\longrightarrow H^{n}_{3-D_\lambda}(\mathfrak g,
V)\longrightarrow H^{n}_{3-Lie}(\mathfrak g, V)\longrightarrow
H_{DO_\lambda}^{n}(d_{\mathfrak
g},d_{V})\longrightarrow\cdot\cdot\cdot.$$

\subsection{Relation between the cohomologies of differential 3-Lie
algebras and differential Leibniz algebras} \hsp We have constructed
cohomological groups of $(\mathfrak g,d_{\mathfrak g})$ with
coefficients in its representation $(V,\rho,d_{V})$. Since any
differential $3$-Lie algebra $(\mathfrak g,d_{\mathfrak g})$
associates with a differential Leibniz algebra
$(L=\wedge^{2}\mathfrak g,d_{L})$. It is natural to investigate the
relation between the cohmologies of differential 3-Lie algebras
$(\mathfrak g,d_{\mathfrak g})$ and associated differential Leibniz
algebra $(L=\wedge^{2}\mathfrak g,d_{L})$. We try to characterize
the relation between them like the case of algebras, unfortunately,
we fail to do that. We only discuss their relationship in the case
of $\lambda=0$.

In the remaining part of this section, we always put $\lambda=0$.

Let $(\mathfrak g,d_{\mathfrak g})$ be a differential 3-Lie algebra
 with weight $\lambda=0$ and $(\rho,V,d_{V})$ be its representation.
 Then $(L=\wedge^{2}\mathfrak g,d_{L})$ is a differential Leibniz
 algebra with weight $\lambda=0$. Furthermore, given linear maps
 $\rho^{L},~\rho^{R}:L\longrightarrow \hbox{Hom}(\mathfrak g,V)$ by
$$\rho^{L}(x,y)(f)(z)=\rho(x,y)f(z)-f([x,y,z])\eqno (3.3)$$ and
$$\rho^{R}(x,y)(f)(z)=f([x,y,z])-\rho(x,y)f(z)-\rho(y,z)f(x)-\rho(z,x)f(y),$$
for any $x \wedge y\in L$. Then $(V,\rho^{L},\rho^{R})$ is a
representation of
Leibniz algebra $L$ \cite{6, 10}.\\

 {\bf Proposition 3.5.} Let $(V,\rho, d_{V})$ be a representation of differential 3-Lie
algebra $(\mathfrak g,d_{\mathfrak g})$, and $(L=\wedge^{2}\mathfrak
g,d_{L})$ be the associated differential Leibniz
 algebra with $(\mathfrak g,d_{\mathfrak g})$, where $d_{L}(x\wedge y)=d_{g}(x)\wedge y+x\wedge
d_{g}(y)$ for any $x\wedge y\in L$. Then $(\hbox{Hom}(\mathfrak
g,V),\rho^{L},\rho^{R},\psi)$ is a representation of the
differential Leibniz
 algebra $(L,d_{L})$, where $\psi:\hbox{Hom}(\mathfrak g,V)\longrightarrow
\hbox{Hom}(\mathfrak g,V)$ is given by
$$\psi(f)=d_{V}f-fd_{\mathfrak g }.\eqno (3.4)$$

\begin{proof} We only need to check that
$$\psi\rho^{L}(X)=\rho^{L}(d_{L}(X))+\rho^{L}(X)\psi$$
and
$$\psi\rho^{R}(X)=\rho^{R}(d_{L}(X))+\rho^{R}(X)\psi.$$
In fact, for any $f\in \hbox{Hom}(\mathfrak g,V)$ and $X=x\wedge
y\in L,~z\in \mathfrak g$, in view of (3.3) and (3.4),
\begin{eqnarray*}
&&\rho^{L}(d_{L}(X))f(z)+\rho^{L}(X)\psi
(f)(z)-(\psi\rho^{L}(X)f)(z)
\\&=&\rho(d_{\mathfrak g}x,y)f(z)-f([d_{\mathfrak g}x,y,z])
+\rho(x,d_{\mathfrak g}x)f(z)-f([x,d_{\mathfrak g}y,z])
+\rho(x,y)d_{V}f(z)\\&&-d_{V}f([x,y,z]) -\rho(x,y)fd_{\mathfrak
g}(z)+f([x,y,d_{\mathfrak g}z])
-d_V(\rho(x,y)f(z)-f([x,y,z]))\\&&+\rho(x,y)fd_{\mathfrak
g}(z)-fd_{\mathfrak g}([x,y,z])\\&=&0.
\end{eqnarray*}
It follows that $$\psi\rho^{L}(X)=\rho^{L}(d_{L
}(X))+\rho^{L}(X)\psi.$$ Analogously,
$$\psi\rho^{R}(X)=\rho^{R}(d_{L}(X))+\rho^{R}(X)\psi.$$
\end{proof}

It is natural to discuss cohomologies of differential Leibniz
algebra $(L ,d_{L})$ with coefficients in the representation
$(\hbox{Hom}(\mathfrak g,V),\rho^{L},\rho^{R},\psi)$.

Let $C_{\hbox{Leib}}^{n}(L,\hbox{Hom}(\mathfrak
g,V))=\hbox{Hom}(\wedge ^{ n} L,\hbox{Hom}(\mathfrak g,V))~~(n\geq
0)$ be the $n$-th cochain group of Leibniz algebra $L$ with
coefficients in the representation $(\hbox{Hom}(\mathfrak
g,V),\rho^{L},\rho^{R},\psi)$ and the coboundary operator
$$\partial_{L}:C_{\hbox{Leib}}^{n}(L,\hbox{Hom}(\mathfrak
g,V))\longrightarrow C_{\hbox{Leib}}^{n+1}(L,\hbox{Hom}(\mathfrak
g,V))$$ is given by
\begin{eqnarray*}
&&\partial_{Leib}(f)(X_1,\cdot\cdot\cdot,X_{n+1})
\\&=&\rho^{L}(X_1)f(X_2,\cdot\cdot\cdot,X_{n+1})+\sum_{i=2}^{n+1}\rho^{R}(X_i)f(X_1,\cdot\cdot\cdot,X_{i-1},X_{i+1},\cdot\cdot\cdot,X_{n+1})
\\&&+\sum_{i< j}(-1)^{j+1}f(X_1,,\cdot\cdot\cdot,X_{i-1},[X_{i},X_{j}]_{F},X_{i+1},\cdot\cdot\cdot,X_{j-1},X_{j+1},\cdot\cdot\cdot,X_{n+1}).
\end{eqnarray*}
for any $X_1,\cdot\cdot\cdot,X_{n+1}\in L$ and $f\in
C_{\hbox{Leib}}^{n}(L,\hbox{Hom}(\mathfrak g,V))$.

Now we are ready to give the cohomology of differential Leibniz
algebra $(L ,d_{L})$ with coefficients in the representation
$(\hbox{Hom}(\mathfrak g,V),\rho^{L},\rho^{R},\psi)$. Denote the
$n$-th cochain group by $C_{DL}^{n}(L,\hbox{Hom}(\mathfrak g,V))$,
where
$$C_{DL}^{n}(L,\hbox{Hom}(\mathfrak
g,V))=C_{\hbox{Leib}}^{n}(L,\hbox{Hom}(\mathfrak g,V))\times
C_{\hbox{Leib}}^{n-1}(L,\hbox{Hom}(\mathfrak g,V))~~(n\geq2), $$ and
$$C_{DL}^{1}(L,\hbox{Hom}(\mathfrak
g,V))=C_{\hbox{Leib}}^{1}(L,\hbox{Hom}(\mathfrak g,V)).$$

Given a linear map
$\delta_{L}:C_{\hbox{Leib}}^{n}(L,\hbox{Hom}(\mathfrak
g,V))\longrightarrow C_{\hbox{Leib}}^{n}(L,\hbox{Hom}(\mathfrak
g,V))$ by $$\delta_{L} f(X_1,\cdot\cdot\cdot, X_{n})
=\sum_{k=1}^{n}f(X_1,\cdot\cdot\cdot, X_{k-1},d_{L}(X_{k}),
X_{k+1},\cdot\cdot\cdot, X_{n})-\psi f(X_1,\cdot\cdot\cdot,
X_{n}).$$ Define the coboundary operator
$$\partial_{DL}:C_{DL}^{n}(L,\hbox{Hom}(\mathfrak g,V))\longrightarrow C_{DL}^{n+1}(L,\hbox{Hom}(\mathfrak
g,V))$$  by
$\partial_{DL}(f,g)=(\partial_{Leib}(f),\partial_{Leib}(g)+(-1)^{n}\delta_{L}
f)$. Then $(C_{DL}^{*}(L,\hbox{Hom}(\mathfrak g,V)),\partial_{DL})$
is a cochain complex, and $\mathcal{H}^{*}_{DL}(L,
\hbox{Hom}(\mathfrak g,V))$ denotes its cohomology group. For more
details on cohomology of Leibniz algebras and differential Leibniz
algebras, one can turn to \cite{6, 7, 10}.

Let $(C_{3-D_\lambda}^{p}(\mathfrak g,V),\partial_{D_\lambda})$ be
the cochain complex of differential $3$-Lie algebra $(\mathfrak
g,d_{\mathfrak g})$ with coefficients in representation
$(V,\rho,d_V)$ studied in Subsection 3.1, and
$\mathcal{H}^{n}_{D_\lambda}(\mathfrak g,V)$ be the corresponding
n-th cohomology group.

 Define $\Theta:C_{3-{\hbox{Lie}}}^{n}(\mathfrak
g,V)\to C_{\hbox{Leib}}^{n-1}(L,\hbox{Hom}(\mathfrak g,V))$ by
$$\Theta(f)(X_1,\cdot\cdot\cdot, X_{n-1})(z)= f(X_1,
\cdot\cdot\cdot,X_{n-1},z)$$ for any $f\in
C_{3-{\hbox{Lie}}}^{n}(\mathfrak g,V)$ and $X_i=x_{i}\wedge y_{i}\in
L$.

 With the above notations, we have the
following result.\\

{\bf Theorem 3.6.} The linear map
$\bar{\Theta}:C_{3-D_\lambda}^{p}(\mathfrak g,V)\longrightarrow
C_{D_\lambda}^{p}(L,\hbox{Hom}(\mathfrak g,V))$ is a cochain map
with $\bar{\Theta}(f,g)=(\Theta f,\Theta g)$ for any $(f,g)\in
C_{3-D_\lambda}^{p}(\mathfrak g,V)$, that is, the following diagram
is commutative:
$$\xymatrix{
 C_{3-D_\lambda}^{n}(\mathfrak g,V) \ar[d]_{\partial_{D_{\lambda}}}
  \ar[r]^{\bar{\Theta}} & C_{DL}^{n-1}(L,\hbox{Hom}(\mathfrak g,V)) \ar[d]^{\partial_{DL}} \\
  C_{3-D_\lambda}^{n+1}(\mathfrak g,V)\ar[r]^{\bar{\Theta}} & C_{DL}^{n}(L,\hbox{Hom}(\mathfrak g,V))
  .}$$
Moreover, $\bar{\Theta}$ induces an isomorphism between
  cohomologies:
$$\mathcal{H}^{n}_{3-D_\lambda}(\mathfrak g, V) \cong\mathcal{H}^{n-1}_{DL}(L, \hbox{Hom}(\mathfrak g,V)).$$

 \begin{proof} For any $(f,g)\in C_{3-D_\lambda}^{n}(\mathfrak
g,V)$, we get
$$\partial_{DL} \bar{\Theta}(f,g)=\partial_{DL}
(\Theta f,\Theta g)=(\partial_{Leib}\Theta f,\partial_{Leib}\Theta
g+(-1)^{n-1}\delta_{L}\Theta f)
$$ and
$$\bar{\Theta}\partial_{D_\lambda}(f,g)=(\Theta \partial f,\Theta
\partial_{\lambda}g+(-1)^{n+1} \Theta\delta f).$$ In the light of
\cite{6, 10}, we only need to check that
 $\delta\Theta=\Theta\delta$.

In fact, for any $x_1,\cdot\cdot\cdot,x_{2n+1}\in \mathfrak g$ and
$f\in C_{3-\hbox{Lie}}^{n}(\mathfrak g,V)$,
\begin{eqnarray*}&&\delta_{L}\Theta(f)(x_1,\cdot\cdot\cdot,x_{2n})(x_{2n+1})
\\&=&\sum_{i=1}^{2n}\Theta(f)(x_1,\cdot\cdot\cdot,x_{i-1},d_{\mathfrak
g}(x_i),x_{i+1},\cdot\cdot\cdot,x_{2n})(x_{2n+1})
-\psi\Theta(f)(x_1,\cdot\cdot\cdot,x_{2n})(x_{2n+1})
\\&=&\sum_{i=1}^{2n}f(x_1,\cdot\cdot\cdot,x_{i-1},d_{\mathfrak
g}(x_i),x_{i+1},\cdot\cdot\cdot,x_{2n},x_{2n+1})
-d_{V}f(x_1,\cdot\cdot\cdot,x_{2n},x_{2n+1})
\\&&+f(x_1,\cdot\cdot\cdot,x_{2n},d_{\mathfrak
g}(x_{2n+1})).
\end{eqnarray*}
On the other hand,
\begin{eqnarray*}&&\Theta\delta(f)(x_1,\cdot\cdot\cdot,x_{2n})(x_{2n+1})
\\&=&\delta(f)(x_1,\cdot\cdot\cdot,x_{2n},x_{2n+1})
\\&=&\sum_{i=1}^{2n+1}f(x_1,\cdot\cdot\cdot,x_{i-1},d_{\mathfrak
g}(x_i),x_{i+1},\cdot\cdot\cdot,x_{2n},x_{2n+1})-d_{V}f(x_1,\cdot\cdot\cdot,x_{2n},x_{2n+1}).
\end{eqnarray*}
Thus, $\Theta\delta=\delta_{L}\Theta$, which yields that
$\bar{\Theta}\partial_{D_\lambda}=\partial_{DL}\bar{\Theta}.$
Clearly, $\bar{\Theta}$ is an isomorphism, the results hold.
 \end{proof}

\section{Deformations of differential 3-Lie algebras}
\hsp In this section, we investigate infinitesimal deformations of a
differential 3-Lie algebra with any weight $\lambda$. We also study
the notion of Nijenhuis operators and $\mathcal{O}$-operators for
differential 3-Lie algebras.

Let $(\mathfrak g,\pi=[\ , \ , \ ],d_{\mathfrak g})$ be a
differential 3-Lie algebra with any weight $\lambda$,
 $\pi_{i}:\wedge^{3}\mathfrak g
 \longrightarrow \mathfrak g
 $ be a trilinear map and $\varphi_{i}:\mathfrak g\longrightarrow
\mathfrak g$ be a linear map.

Consider the space $\mathfrak g[[t]]$ of formal power series in $t$
with coefficients in $\mathfrak g$ and a $t$-parametrized family of
trilinear operations
 $$\pi_{t}(x,y,z)=\sum_{i=0}^{n-1}t^{i}\pi_{i}(x,y,z),$$
 and
linear operations
 $$\varphi_{t}(x)=\sum_{i=0}^{n-1}t^i\varphi_{i}(x),$$
 where $\varphi_{0}=d_{\mathfrak g}$ and $\pi_{0}=\pi$.

 If all $(\mathfrak g[[t]],
 \pi_{t},\varphi_{t})$ are differential 3-Lie algebras, we say that
 $(\pi_{t},\varphi_t)$ generates a deformation of the differential 3-Lie algebra $(
 \mathfrak g,\pi,d_{\mathfrak g})$.

If all $(\mathfrak g[[t]],
 \pi_{t},\varphi_{t})$ are differential 3-Lie
algebras with weight $\lambda$, then we have
\begin{eqnarray*}&&\sum_{i+j=k,~i,~j\geq
 0}(\pi_i(\pi_j(v,w,x),y,z)+\pi_i(x,\pi_j(v,w,y),z)+\pi_i(x,y,\pi_j(v,w,z))\\&&-\pi_i(v,w,\pi_j(x,y,z)))=0, ~~~~~~~~~
 ~~~~~~~~~~~~~~~~~~~~~~~~~~~~~~~~~~~~~~~~~~~~~~~~~~~~~~~~~~~~~(4.1)
\end{eqnarray*} and
\begin{eqnarray*}&&\sum_{i+j=n,~i,~j\geq
 0}(\pi_i(\varphi_{j}x,y,z)+\pi_i(x,\varphi_{j}y,z)+\pi_i(x,y,\varphi_{j}z))
\\&&+\lambda\sum_{i+j+k=n,~i,~j,~k\geq
 0}(\pi_i(\varphi_{j}x,\varphi_{k}y,z)+\pi_i(x,\varphi_{j}y,\varphi_{k}z)+\pi_i(\varphi_{j}x,y,\varphi_{k}z))
\\&&+\sum_{i+j+k+l=n,~i,~j,~k,~l\geq
 0}\lambda^{2}\pi_i(\varphi_{j}x,\varphi_{k}y,\varphi_{l}z)-\sum_{i+j=n,~i,~j\geq
 0}\varphi_{i}\pi_j(x,y,z)=0
.~~~~~~~(4.2)\end{eqnarray*}\\

 {\bf Proposition 4.1.}
 $(\mathfrak g[[t]],
 \pi_{t},\varphi_{t})$ generates a deformation of the differential 3-Lie algebra $(
 \mathfrak g,\pi,d_{\mathfrak g})$. Then $(\pi_1,\varphi_1)$ is a 2-cocycle of $(
 \mathfrak g,d_{\mathfrak g})$ with coefficients in the adjoint
 representation.

\begin{proof} If $(\mathfrak g[[t]],
 \pi_{t},\varphi_{t})$ generates a deformation of the differential 3-Lie algebra $(
 \mathfrak g,\pi,d_{\mathfrak g})$, then $(4.1)$ and $(4.2)$ hold.

On the other hand, for any $(\pi_1,\varphi_1)\in
C_{3-D_\lambda}^{2}(\mathfrak g,\mathfrak
g)=C_{3-{\hbox{Lie}}}^{2}(\mathfrak g,\mathfrak g)\times
C_{3-{\hbox{DLie}}}^{1}(d_\mathfrak g,d_\mathfrak g)$ and $X=x\wedge
y\in L$, $z\in \mathfrak g$, $(\pi_1,\varphi_1)$ is a 2-cocycle if
$$\partial_{D_\lambda}(\pi_1,\varphi_1)=(\partial \pi_1,\partial_\lambda
\varphi_1+(-1)^{2}\delta \pi_1)=0,$$ that is, $\partial \pi_1=0$ and
$\partial_\lambda \varphi_1+\delta \pi_1=0$. From the deformation
theory of 3-Lie algebras \cite{34}, we know that $\partial \pi_1=0$
is equivalent to (4.1) when $k=1$. By computation,
\begin{eqnarray*}
&&\partial_\lambda
(\varphi_1)(x,y,z)\\&=&-\pi_1([x,y,z])+\hat{\rho}(x,y)\varphi_1(z)+\hat{\rho}(y,z)\varphi_1(x)+\hat{\rho}(z,x)\varphi_1(y)\\
&=&-\varphi_1([x,y,z])+[x,y,\varphi_1(z)]+\lambda[d_{\mathfrak
g}x,y,\varphi_1(z)]+\lambda[x,d_{\mathfrak
g}y,\varphi_1(z)]+\lambda^{2}[d_{\mathfrak g}x,d_{\mathfrak
g}y,\varphi_1(z)]\\&&+ [y,z,\varphi_1(x)]+\lambda[d_{\mathfrak
g}y,z,\varphi_1(x)]+\lambda[y,d_{\mathfrak
g}z,\varphi_1(x)]+\lambda^{2}[d_{\mathfrak g}y,d_{\mathfrak
g}z,\varphi_1(x)]\\&&+ [z,x,\varphi_1(y)]+\lambda[d_{\mathfrak
g}z,x,\varphi_1(y)]+\lambda[z,d_{\mathfrak
g}x,\varphi_1(y)]+\lambda^{2}[d_{\mathfrak g}z,d_{\mathfrak
g}x,\varphi_1(y)]
\end{eqnarray*}
and
\begin{eqnarray*}
&&\delta \pi_1(x,y,z)\\&=& \pi_1(d_{\mathfrak
g}x,y,z)+\pi_1(x,d_{\mathfrak g}y,z)+\pi_1(x,y,d_{\mathfrak
g}z)+\lambda(\pi_1(d_{\mathfrak g}x,d_{\mathfrak
g}y,z)\\&&+\pi_1(x,d_{\mathfrak g}y,d_{\mathfrak
g}z)+\pi_1(d_{\mathfrak g}x,y,d_{\mathfrak
g}z))+\lambda^2\pi_1(d_{\mathfrak g}x,d_{\mathfrak g}y,d_{\mathfrak
g}z)-d_{\mathfrak g}\pi_1(x,y,z).
\end{eqnarray*}
These imply that $\partial_\lambda \varphi_1+\delta \pi_1=0$ is
equivalent to (4.2) when $n=1$.
\end{proof}

{\bf Definition 4.2.} A
 deformation of a differential 3-Lie algebra $(\mathfrak g,d_{\mathfrak g
 })$ is said to be trivial if there is a linear map $N:\mathfrak g
 \longrightarrow \mathfrak g
 $ such that $K_{t}=I+tN~(\forall~ t)$ satisfies
 $$K_{t}d_{\mathfrak g
 }=d_{\mathfrak g
 }K_{t}$$ and
 $$K_{t}[x,y,z]_{t}=[K_{t}x, K_{t}y,K_{t}z]
 .$$

Nijenhuis operator and $\mathcal{O}$-operator on 3-Lie algebras are
discussed in \cite{14, 33}. On the basis of these works, we will
consider theses operators on differential 3-Lie algebras.\\

 {\bf Definition 4.3.} Let $(\mathfrak g,d_{\mathfrak g
 })$ be a differential 3-Lie algebra with weight $\lambda$. A linear map $N:\mathfrak g
 \longrightarrow \mathfrak g $ is called a Nijenhuis operator on the differential 3-Lie algebra $(\mathfrak g,d_{\mathfrak g
 })$
 if the following hold:
 $$Nd_{\mathfrak g
 }=d_{\mathfrak g
 }N$$
 and
\begin{eqnarray*}[Nx,Ny,Nz]
 &=&N([Nx,Ny,z]+[x,Ny,Nz]+[Nx,y,Nz])\\&&-N^{2}([Nx,y,z]+[x,Ny,z]+[x,y,Nz])+N^{3}([x,y,z]).
 \end{eqnarray*}

By direct calculation, we have\\

 {\bf Proposition 4.4.} If $(\mathfrak g,d_{\mathfrak g
 })$ is a differential 3-Lie algebra, then
$(\mathfrak g,[ \ , \ , \ ]_{N},d_{\mathfrak g
 })$ is also a differential 3-Lie
algebra, where
\begin{eqnarray*}[x,y,z]_{N}&=&[Nx,Ny,z]+[x,Ny,Nz]+[Nx,y,Nz]\\&&-N([Nx,y,z]+[x,Ny,z]+[x,y,Nz])+N^{2}([x,y,z]).
\end{eqnarray*}

{\bf Definition 4.5.} A linear map $K:V\longrightarrow \mathfrak g$
 is called an $\mathcal{O}$-operator on the differential 3-Lie algebra $(\mathfrak g,d_{\mathfrak g})$
  associated to the representation $(V,\rho,d_{V})$ if $K$ is an
$\mathcal{O}$-operator on the 3-Lie algebra $\mathfrak g$ associated
to the representation $(V,\rho)$, that is, for all $u,v,w\in V$,
$$[K(u), K(v), K(w)]=K(\rho (K(u),K(v))w+\rho(K(v),K(w))u+\rho (K(w),K(u))v)\eqno (4.3)$$
 and $$Kd_{V}
 =d_{\mathfrak g}K.\eqno (4.4)$$

{\bf Proposition 4.6.} Let $K:V\longrightarrow \mathfrak g$ be an
$\mathcal{O}$-operator on the differential 3-Lie algebra $(\mathfrak
g,d_{\mathfrak g})$ associated to the representation
$(V,\rho,d_{V})$. Then $(V,[ \ , \ , \ ]_{K},d_{V})$ is a
differential 3-Lie algebra, where
 $$ [u,v,w]_{K}=\rho (K(u),K(v))w+\rho(K(v),K(w))u+\rho (K(w),K(u))v \eqno (4.5)$$ for all $u,v,w \in
 V$.
\begin{proof} It can be obtained by direct calculation.
\end{proof}

{\bf Proposition 4.7.} $K:V\longrightarrow \mathfrak g$ is an
$\mathcal{O}$-operator on the differential 3-Lie algebra $(\mathfrak
g,d_{\mathfrak g})$ associated to the representation
$(V,\hat{\rho},d_{V})$ if and only if  $\hat{K}=K+\lambda
d_{\mathfrak g} K$ is an $\mathcal{O}$-operator on the differential
3-Lie algebra $(\mathfrak g,d_{\mathfrak g})$ associated to the
representation $(V,\rho,d_{V})$, where $\hat{\rho}$ is mentioned in
Section 2.

\begin{proof}
Notice that
\begin{eqnarray*}&&\hat{\rho} (K(u),K(v))\\&=&\rho(K(u),K(v))+\lambda(\rho(d_{\mathfrak
g}K(u),K(v))+\lambda\rho(K(u),d_{\mathfrak
g}K(v))+\lambda^{2}\rho(d_{\mathfrak g}K(u),d_{\mathfrak
g}K(v))\\&=& \rho(\hat{K}(u),\hat{K}(v)).
\end{eqnarray*}
So we get the result.
\end{proof}

{\bf Proposition 4.8.} Let $K:V\longrightarrow \mathfrak g$ be an
$\mathcal{O}$-operator on the differential 3-Lie algebra $(\mathfrak
g,d_{\mathfrak g})$ associated to the representation
$(V,\rho,d_{V})$. Define
 linear maps
$\varrho_{K}:V\wedge V\longrightarrow {gl}(\mathfrak g)$ by
$$\varrho_{K}(u,v)x=[K(u),K(v),x]-K(\rho(Kv,x)u+\rho(x,Ku)v).$$
Then $(\mathfrak g,\varrho_{K},d_{\mathfrak g})$ is a representation
of $(V,[ \ , \ , \  ]_{K},d_{V})$.

\begin{proof} In view of \cite{33}, $(\mathfrak g,\varrho_{K})$ is a representation of $(V,[
\ , \ , \  ]_{K})$. We only need to check that (2.5) holds for
$\varrho_{K}$.
\begin{eqnarray*}
&&(\varrho_{K}(d_{V}u,v)+\varrho_{K}(u,d_{V}v)+\lambda\varrho_{K}(d_{V}u,d_{V}v)+\varrho_{K}(u,v)d_{\mathfrak
g})x+\lambda(\varrho_{K}(d_{V}u,v)
\\&&+\varrho_{K}(u,d_{V}v)+\lambda\varrho_{K}(d_{V}u,d_{V}v))d_{\mathfrak
g}x
\\
&=&[Kd_{V}u,Kv,x]-K(\rho(Kv,x)d_{V}u+\rho(x,Kd_{V}u)v)+[Ku,Kd_{V}v,x]-K(\rho(Kd_{V}v,x)u\\&&+\rho(x,Ku)d_{V}v)
+\lambda([Kd_{V}u,Kd_{V}v,x]-K(\rho(Kd_{V}v,x)d_{V}u+\rho(x,Kd_{V}u)d_{V}v))
\\&&+[Ku,Kv,d_{\mathfrak
g}x]-K(\rho(Kv,d_{\mathfrak g}x)u+\rho(d_{\mathfrak g}x,Ku)v)
+\lambda([Kd_{V}u,Kv,d_{\mathfrak g}x]\\&&-K(\rho(Kv,d_{\mathfrak
g}x)d_{V}u+\rho(d_{\mathfrak g}x,Kd_{V}u)v)
+[Ku,Kd_{V}v,d_{\mathfrak g}x]-K(\rho(Kd_{V}v,d_{\mathfrak
g}x)u\\&&+\rho(d_{\mathfrak g}x,Ku)d_{V}v))
+\lambda^2([Kd_{V}u,Kd_{V}v,d_{\mathfrak
g}x]-K(\rho(Kd_{V}v,d_{\mathfrak g}x)d_{V}u+\rho(d_{\mathfrak
g}x,Kd_{V}u)d_{V}v))
\\
&=&[d_{\mathfrak g}Ku,Kv,x]-K(\rho(Kv,x)d_{V}u+\rho(x,d_{\mathfrak
g}Ku)v)+[Ku,d_{\mathfrak g}Kv,x]-K(\rho(d_{\mathfrak
g}Kv,x)u\\&&+\rho(x,Ku)d_{V}v) +\lambda([d_{\mathfrak
g}Ku,d_{\mathfrak g}Kv,x]-K(\rho(d_{\mathfrak
g}Kv,x)d_{V}u+\rho(x,d_{\mathfrak g}Ku)d_{V}v))
\\&&+[Ku,Kv,d_{\mathfrak
g}x]-K(\rho(Kv,d_{\mathfrak g}x)u+\rho(d_{\mathfrak g}x,Ku)v)
+\lambda([d_{\mathfrak g}Ku,Kv,d_{\mathfrak
g}x]\\&&-K(\rho(Kv,d_{\mathfrak g}x)d_{V}u+\rho(d_{\mathfrak
g}x,d_{\mathfrak g}Ku)v) +[Ku,d_{\mathfrak g}Kv,d_{\mathfrak
g}x]-K(\rho(d_{\mathfrak g}Kv,d_{\mathfrak
g}x)u\\&&+\rho(d_{\mathfrak g}x,Ku)d_{V}v)) +\lambda^2([d_{\mathfrak
g}Ku,d_{\mathfrak g}Kv,d_{\mathfrak g}x]-K(\rho(d_{\mathfrak
g}Kv,d_{\mathfrak g}x)d_{V}u+\rho(d_{\mathfrak g}x,d_{\mathfrak
g}Ku)d_{V}v))
\\
&=& d_{\mathfrak g}[Ku,Kv,x]-d_{\mathfrak g}K
(\rho(Kv,x)u+\rho(x,Ku)v)
\\
&=& d_{\mathfrak g}(\varrho_{K}(u,v)x).
\end{eqnarray*}
\end{proof}
Consider the cohomology of the new differential 3-Lie algebra $(V,[
\ , \ , \  ]_{K},d_{V})$ with coefficients in the representation
$(\mathfrak g,\varrho_{K},d_{\mathfrak g})$.

Denote the associated Leibniz algebra by $(L(V)=\wedge^{2}V, [ \ , \
]_{VF})$ with the Leibniz bracket $[ \ , \  ]_{VF}$ given by
$$[ U ,V]_{VF}=v_1\wedge [u_1,u_2,v_2]+[u_1,u_2,v_1]\wedge v_2
$$ for all $U=u_1\wedge u_2\in L(V) $ and
$V=v_1\wedge v_2\in L(V)$. The space
$C_{3-{\hbox{Lie}}}^{p}(V,\mathfrak g)$ of $p$-cochains ($p\geq 1$)
is the set of multilinear maps of the form
$$f:\wedge^{p-1}L(V)\wedge V\longrightarrow \mathfrak g,$$
and the coboundary operator
$\partial:C_{3-{\hbox{Lie}}}^{p}(V,\mathfrak g)\longrightarrow
C_{3-{\hbox{Lie}}}^{p+1}(V,\mathfrak g)$ is as follows:
\begin{eqnarray*}&&(\partial f)(V_1,\cdot\cdot\cdot,V_{p},w)\\&=&\sum_{1\leq i<
k\leq {p}}
(-1)^{i}f(V_1,\cdot\cdot\cdot,\hat{V_i},\cdot\cdot\cdot,V_{k-1},[V_i,V_k]_{VF},V_{k+1},\cdot\cdot\cdot,V_{p}
,w)
\\&&+\sum_{i=1}^{p}(-1)^{i}f(V_1,\cdot\cdot\cdot,\hat{V_i},\cdot\cdot\cdot,V_{p},[V_i,w]_{K})+
\sum_{i=1}^{p}(-1)^{i+1}\varrho_{K}(V_i)f(V_1,\cdot\cdot\cdot,\hat{V_i},\cdot\cdot\cdot,V_{p},w)\\&&
+(-1)^{p+1}(\varrho_{K}(v_{p},w)f(V_1,\cdot\cdot\cdot,V_{p-1},u_{p})+\varrho_{K}(w,u_{p})f(V_1,\cdot\cdot\cdot,V_{p-1},v_{p}))
\end{eqnarray*}
for all $V_i=u_{i}\wedge v_{i}\in V\wedge V , w\in V$.

Moreover, let $(\mathfrak g,\hat{\varrho_{K}},d_{\mathfrak g})$
denote the new representation defined as in Proposition 2.6.

Denote $C_{3-{\hbox{DLie}}}^{p}(d_V,d_\mathfrak
g)=\hbox{Hom}(\wedge^{p-1}L(V)\wedge V,\mathfrak g)$ by the set of
$p$-cochains ($p\geq 1$) of the differential operator $d_V$ with
coefficients in the representation $(\mathfrak
g,\hat{\varrho_{K}},d_\mathfrak g)$, and the coboundary operator
$$\partial_\lambda:C_{3-{\hbox{DLie}}}^{p}(d_V,d_\mathfrak
g)\longrightarrow C_{3-{\hbox{DLie}}}^{p+1}(d_V,d_\mathfrak g)$$ is
given by
\begin{eqnarray*}&&(\partial_\lambda f)(V_1,\cdot\cdot\cdot,V_{p},w)\\&=&
\sum_{1\leq i< k\leq {p}}
(-1)^{i}f(V_1,\cdot\cdot\cdot,\hat{V_i},\cdot\cdot\cdot,V_{k-1},[V_i,V_k]_{VF},V_{k+1},\cdot\cdot\cdot,V_{p}
,w)
\\&&+\sum_{i=1}^{p}(-1)^{i}f(V_1,\cdot\cdot\cdot,\hat{V_i},\cdot\cdot\cdot,V_{p},[V_i,w]_{K})+
\sum_{i=1}^{p}(-1)^{i+1}\hat{\varrho_{K}}(V_i)f(V_1,\cdot\cdot\cdot,\hat{V_i},\cdot\cdot\cdot,V_{p},w)\\&&
+(-1)^{p+1}(\hat{\varrho_{K}}(v_{p},w)f(V_1,\cdot\cdot\cdot,V_{p-1},u_{p})+\hat{\varrho_{K}}(w,u_{p})f(V_1,\cdot\cdot\cdot,V_{p-1},v_{p}))
\end{eqnarray*}
for any $f\in C_{3-{\hbox{DLie}}}^{p}(d_V,d_\mathfrak g)$,
$V_1,\cdot\cdot\cdot,V_{p}\in V\wedge V,~w\in V$.\\

{\bf Definition 4.9.} The cohomology of the cochain complex
$(C_{3-{\hbox{DLie}}}^{*}(d_V,d_\mathfrak g),\partial_\lambda)$,
denoted by $H_{DO_\lambda}^{*}(d_V,d_{\mathfrak g})$, is called the
cohomology of the differential operator $d_V$ with coefficients in
the representation $(\mathfrak
g_{\lambda},\hat{\varrho_{K}},d_\mathfrak g)$.
\\

Given a linear map $\delta:C_{3-{\hbox{Lie}}}^{p}(V,\mathfrak
g)\longrightarrow C_{3-{\hbox{DLie}}}^{p}(d_V,d_\mathfrak g)$ by
$$\delta f(V_1,\cdot\cdot\cdot,V_{p},w)=
\sum_{k=1}^{2p+1}\lambda^{k-1}f^{k}(V_1,
\cdot\cdot\cdot,V_{p},w)-d_{\mathfrak g}f(V_1,
\cdot\cdot\cdot,V_{p},w).
$$
Denote the set of $p$-cochains by
$$C_{3-D_\lambda}^{p}(V,\mathfrak g)=:C_{3-{\hbox{Lie}}}^{p}(V,\mathfrak g)\times
C_{3-{\hbox{DLie}}}^{p-1}(d_V,d_\mathfrak g),~~p\geq 2$$ and
$$C_{3-D_\lambda}^{1}(V,\mathfrak g)=\hbox{Hom}(V,\mathfrak g).$$

Define the coboundary operator
$\partial_{D_\lambda}:C_{3-D_\lambda}^{p}(V,\mathfrak
g)\longrightarrow C_{3-D_\lambda}^{p+1}(V,\mathfrak g)$ by
$$\partial_{D_\lambda}(f,g)=(\partial f,\partial_\lambda
g+(-1)^{p}\delta f)$$ for any $(f,g)\in
C_{3-D_\lambda}^{p}(V,\mathfrak g)$.

Associated to the representation $(\mathfrak g,d_{\mathfrak g})$, we
obtain a cochain complex $(C_{3-D_\lambda}^{*}(V,\mathfrak g),
\partial_{D_\lambda})$. Denote the cohomology group of this cochain
complex by $\mathcal{H}^{*}_{3-D_\lambda}(V,\mathfrak g)$, which is
called the cohomology of the differential $3$-Lie algebra
$(V,d_{V})$ with coefficients in representation $(\mathfrak
g,\varrho_{K},d_{\mathfrak g})$.

We calculate the 1-cocycle.

For any $f\in \hbox{Hom}(V,\mathfrak g)$,
$$\partial_{D_\lambda}(f)=(\partial f,-\delta f),$$
where
\begin{eqnarray*}&&\partial
(f)(u,v,w)\\&=&-f([u,v,w]_{K})+\varrho_{K}(u,v)f(w)+\varrho_{K}(v,w)f(u)+\varrho_{K}(w,u)f(v)
\\&=&-f(\rho (K(u),K(v))w+\rho(K(v),K(w))u+\rho (K(w),K(u))v )+[Ku,Kv,f(w)]\\&&-K(\rho(Kv,f(w))u+\rho(f(w),Ku)v
+[Kv,Kw,f(u)]-K(\rho(Kw,f(u))v\\&&+\rho(f(u),Kv)w
+[Kw,Ku,f(v)]-K(\rho(Ku,f(v))w+\rho(f(v),Kw)u,
\end{eqnarray*}
and $$\delta f(w)=fd_{V}(w)-d_\mathfrak g f(w)$$ for any $u\wedge
v\in V\wedge V,w\in V $.

It follows that $f\in \hbox{Hom}(V,\mathfrak g)$ is a 1-cocycle if
and only if
\begin{eqnarray*}&&f(\rho (K(u),K(v))w+\rho(K(v),K(w))u+\rho (K(w),K(u))v )\\&=&[Ku,Kv,f(w)]-K(\rho(Kv,f(w))u+\rho(f(w),Ku)v
+[Kv,Kw,f(u)]\\&&-K(\rho(Kw,f(u))v+\rho(f(u),Kv)w
+[Kw,Ku,f(v)]-K(\rho(Ku,f(v))w\\&&+\rho(f(v),Kw)u ;
\end{eqnarray*}
 and $$fd_{V}=d_\mathfrak g f.$$

Thus, we have\\

{\bf Proposition 4.10.} The $\mathcal{O}$-operator
$K:V\longrightarrow \mathfrak g$ on differential $3$-Lie algebra
$(\mathfrak g,d_{\mathfrak g})$
 associated to the representation $(V,\rho,d_{V})$ is a
 1-cocycle of the cochain complex $(C_{3-D_\lambda}^{*}(V,\mathfrak
g),\partial_{D_\lambda})$.

\section{Abelian extensions of differential 3-Lie algebras}
\hsp Let $(\mathfrak g,d_{\mathfrak g})$ be a differential 3-Lie
algebra with any weight $\lambda$ and $(V,d_{V})$ an
abelian differential 3-Lie algebra with the trivial product.\\

{\bf Definition 5.1.} A central extension of $(\mathfrak
g,d_{\mathfrak g})$ by the abelian differential 3-Lie algebra
$(V,d_{V})$ is an exact sequence of differential 3-Lie algebras
$$0\longrightarrow(V,d_{V})\stackrel{i}{\longrightarrow} (\hat{\mathfrak g},d_{\hat{\mathfrak g}})
\stackrel{p}{\longrightarrow} (\mathfrak g,d_{\mathfrak
g})\longrightarrow0$$ such that $[V,V,\hat{\mathfrak
g}]=[V,\hat{\mathfrak g},V]=[\hat{\mathfrak g},V,V]=0$.\\

{\bf Definition 5.2.} Let $(\hat{{\mathfrak g}}_1,d_{\hat{{\mathfrak
g}}_1})$ and $(\hat{{\mathfrak g}}_2,d_{\hat{{\mathfrak g}}_2})$ be
two abelian extensions of $(\mathfrak g,d_{\mathfrak g})$ by
$(V,d_{V})$. They are said to be equivalent if there is a
homomorphism of differential 3-Lie algebras
$\varphi:(\hat{{\mathfrak g}}_1,d_{\hat{{\mathfrak
g}}_1})\longrightarrow (\hat{{\mathfrak g}}_2,d_{\hat{{\mathfrak
g}}_2})$ such that the following commutative diagram holds:
$$\xymatrix{
  0 \ar[r] & (V,d_{V}) \ar[d]_{id} \ar[r]^{i} & (\hat{{\mathfrak g}}_1,d_{\hat{{\mathfrak g}}_1})
  \ar[d]_{\varphi} \ar[r]^{p} & ({\mathfrak g},d_{{\mathfrak g}}) \ar[d]^{id} \ar[r] & 0\\
 0 \ar[r] & (V,d_{V}) \ar[r]^{i} & (\hat{{\mathfrak g}}_2,d_{\hat{{\mathfrak g}}_2}) \ar[r]^{p} & ({\mathfrak g},d_{{\mathfrak g}})  \ar[r] & 0
 .}$$

Suppose $(\hat{\mathfrak g},d_{\hat{\mathfrak g}})$ is an abelian
extension of $(\mathfrak g,d_{\mathfrak g})$ by $(V,d_{V})$. Let
$s:\mathfrak g\longrightarrow \hat{\mathfrak g}$ be any section of
$p$. Define bilinear map $\rho:{\mathfrak g}\wedge {\mathfrak
g}\longrightarrow gl(V)$ by
$$\rho(x,y)v=[s(x),s(y),v]_{\hat{\mathfrak g}},$$
and trilinear map $\psi:\wedge^{3} {\mathfrak g}\longrightarrow V$
and linear map $\chi:{\mathfrak g}\longrightarrow V$ respectively by
$$\psi(x,y,z)=[s(x),s(y),s(z)]-s[x,y,z],$$
$$\chi(x)=d_{\hat{\mathfrak g}}(s(x))-s(d_{\mathfrak g}(x)).$$

{\bf Theorem 5.3.} With the above notations, $(V,\rho,\theta,D,d_V)$
is a representation of $(\mathfrak g,d_{\mathfrak g})$ and it is
independent on the choice of section $s$. Furthermore, equivalent
abelian extensions give the same representation.

\begin{proof} From the case of 3-Lie algebras \cite{34}, we only need to check that (2.5) holds.
For any $x,y,z\in \mathfrak g$, $v\in V$, since $ d_{\hat{\mathfrak
g}}(s(x))-s(d_{\mathfrak g}(x)),~s[x,y,z]-[s(x), s(y), s(z)]\in V$
and $V$ is abelian, which yield that (2.5) holds.
\end{proof}

{\bf Proposition 5.4.} The vector space $(\mathfrak g\oplus V,[ \ ,
\ , \ ]_{\psi},d_{\chi})$ with the bracket
$$[x+a,y+b,z+c]_{\psi}=[x,y,z]+\psi(x,y,z)+\rho(x,y)c+\rho(y,z)a+\rho(z,x)b\eqno (5.1)$$ and
$$d_{\chi}(x+a)=d_{\mathfrak g}(x)+d_{V}(a)+\chi(x)\eqno (5.2)$$ is a
differential 3-Lie algebra with weight $\lambda$ if and only if
$(\psi,\chi)$ is a 2-cocycle of $(\mathfrak g,d_{\mathfrak g})$ with
coefficients in the representation $(V,d_{V})$.

\begin{proof} For any $x,y,z\in \mathfrak g,~a,b,c\in V$, from
the equality
\begin{eqnarray*}&&
d_{\chi}([x+a,y+b,z+c]_{\psi})\\&=&[d_{\chi}(x+a),y+b,z+c]_{\psi}
+[x+a,d_{\chi}(y+b),z+c]_{\psi} +[x+a,y+b,d_{\chi}(z+c)]_{\psi}
\\&&+\lambda([d_{\chi}(x+a),d_{\chi}(y+b),z+c]_{\psi}
+[x+a,d_{\chi}(y+b),d_{\chi}(z+c)]_{\psi} \\&&+
[d_{\chi}(x+a),y+b,d_{\chi}(z+c)]_{\psi} +
\lambda[d_{\chi}(x+a),d_{\chi}(y+b),d_{\chi}(z+c)]_{\psi}),
 \end{eqnarray*}
we have
\begin{eqnarray*}&&\chi([x,y,z])+d_{V}\psi(x,y,z)\\&=&
\psi(d_{\mathfrak g}x,y,z) +\psi(x,d_{\mathfrak
g}y,z)+\psi(x,y,d_{\mathfrak g}z)
 +\lambda(\psi(d_{\mathfrak g}x,d_{\mathfrak g}y,z)
  +\psi(x,d_{\mathfrak g}y,d_{\mathfrak g}z)
\\&&+\psi(d_{\mathfrak g}x,y,d_{\mathfrak g}z) +\lambda\psi(d_{\mathfrak g}x,d_{\mathfrak g
}y,d_{\mathfrak
g}z))+\lambda(\rho(y,z)\chi(x)+\rho(z,x)\chi(y)+\rho(x,y)\chi(z)
\\&&+\rho(d_{\mathfrak g}y,z)\chi(x)+\rho(z,d_{\mathfrak g}x)\chi(y)
+\rho(d_{\mathfrak g}z,x)\chi(y)+\rho(x,d_{\mathfrak g}y)\chi(z)
+\rho(y,d_{\mathfrak g}z)\chi(x)\\&&+\rho(d_{\mathfrak
g}x,y)\chi(z)) +\lambda^{2}(\rho(d_{\mathfrak g}y,d_{\mathfrak
g}z)\chi(x)+\rho(d_{\mathfrak g}x,d_{\mathfrak
g}y)\chi(z)+\rho(d_{\mathfrak g}x,d_{\mathfrak g}z)\chi(y)).~(5.3)
 \end{eqnarray*}
On the other hand, for any $(\psi,\chi)\in
C_{3-D_\lambda}^{2}(\mathfrak g,\mathfrak g)$ is a 2-cocycle if
$$\partial_{D_\lambda}(\psi,\chi)=(\partial \psi,\partial_\lambda
\chi+(-1)^{2}\delta \psi)=0,$$ that is, $\partial \psi=0$ and
$\partial_\lambda \chi+\delta \psi=0$.

Based on the deformation theory of 3-Lie algebras \cite{34}, we only
need to prove that $\partial_\lambda \chi+\delta \psi=0$ if and only
if (5.3) holds. It can be checked directly.
\end{proof}
\section{Skeletal differential 3-Lie 2-algebras and crossed modules}
\hsp In this section, we introduce the notion of weighted
differential 3-Lie 2-algebras and show that skeletal differential
3-Lie 2-algebras are classified by 3-cocycles of weighted
differential 3-Lie algebras.

 We start with recalling the definition of 3-Lie 2-algebras
 \cite{101}.

 A 3-Lie 2-algebra is a quintuple $(\mathfrak
g_0,\mathfrak g_1,h,l_3,l_5)$, where $h:\mathfrak
g_1\longrightarrow\mathfrak g_0$ is a linear map, $l_3:\mathfrak
g_i\wedge \mathfrak g_j\wedge\mathfrak g_k\longrightarrow\mathfrak
g_{i+j+k}~~(0\leq i+j+k\leq 1)$ are completely skew-symmetric
trilinear maps and $l_5:\wedge^{2}\mathfrak g_0 \wedge
\wedge^{3}\mathfrak g_0\longrightarrow\mathfrak g_{1}$ is a
multilinear map, and for any $x_i\in \mathfrak
g_0~(i=1,\cdot\cdot\cdot,7),~a\in \mathfrak
 g_1$, the followings are hold:
$$hl_3(x,y,a)=l_3(x,y,h(a)),~~~l_3(a,b,x)=l_3(h(a),b,x),\eqno(6.1)$$
\begin{eqnarray*}hl_5(x_1,x_2,x_3,x_4,x_5)&=&-l_3(x_1,x_2,l_3(x_3,x_4,x_5))+l_3(x_3,l_3(x_1,x_2,x_4),x_5)
\\&&+l_3(l_3(x_1,x_2,x_3),x_4,x_5)+l_3(x_3,x_4,l_3(x_1,x_2,x_5)),~~~~~~~~~~(6.2)
 \end{eqnarray*}
\begin{eqnarray*}l_5(h(a),x_2,x_3,x_4,x_5)&=&-l_3(a,x_2,l_3(x_3,x_4,x_5))+l_3(x_3,l_3(a,x_2,x_4),x_5)
\\&&+l_3(l_3(a,x_2,x_3),x_4,x_5)+l_3(x_3,x_4,l_3(a,x_2,x_5)),~~~~~~~~~~~~(6.3)
 \end{eqnarray*}
\begin{eqnarray*}l_5(x_1,x_2,h(a),x_4,x_5)&=&-l_3(x_1,x_2,l_3(a,x_4,x_5))+l_3(a,l_3(x_1,x_2,x_4),x_5)
\\&&+l_3(l_3(x_1,x_2,a),x_4,x_5)+l_3(a,x_4,l_3(x_1,x_2,x_5)),~~~~~~~~~~~~(6.4)
 \end{eqnarray*}
\begin{eqnarray*}&&l_3(l_5(x_1,x_2,x_3,x_4,x_5),x_6,x_7)+l_3(x_5,l_5(x_1,x_2,x_3,x_4,x_6),x_7)
\\&&+l_3(x_1,x_2,l_5(x_3,x_4,x_5,x_6,x_7))+l_3(x_5,x_6,l_5(x_1,x_2,x_3,x_4,x_7))
\\&&+l_5(x_1,x_2,l_3(x_3,x_4,x_5),x_6,x_7)+l_5(x_1,x_2,x_5,l_3(x_3,x_4,x_6),x_7)
\\&&+l_5(x_1,x_2,x_5,x_6,l_3(x_3,x_4,x_7))=l_3(x_3,x_4,l_5(x_1,x_2,x_5,x_6,x_7))
\\&&+l_5(l_3(x_1,x_2,x_3),x_4,x_5,x_6,x_7)
+l_5(x_3,l_3(x_1,x_2,x_4),x_5,x_6,x_7)\\&&+l_5(x_3,x_4,l_3(x_1,x_2,x_5),x_6,x_7)
+l_5(x_3,x_4,x_5,x_6,l_3(x_1,x_2,x_7))
\\&&+l_5(x_1,x_2,x_3,x_4,l_3(x_5,x_6,x_7))
+l_5(x_3,x_4,x_5,x_6,l_3(x_1,x_2,x_7)).~~~~~~~~~~~~~~~~~~~~(6.5)
 \end{eqnarray*}
 A 3-Lie 2-algebra is called skeletal (strict) if $h=0$
($l_5=0$).

Motivated by \cite{32} and \cite{102}, we give the definition of a
differential 3-Lie 2-algebra with any weight
$\lambda$. \\

{\bf Definition 6.1.} A differential 3-Lie 2-algebra with any weight
$\lambda$ consists of a 3-Lie 2-algebra $\mathfrak g=(\mathfrak
g_0,\mathfrak g_1,h,l_3,l_5)$
 and a differential operator $d=(d_0,d_1,d_2)$ of $\mathfrak g$,
 where $d_0:\mathfrak g_0\longrightarrow \mathfrak g_0,~~d_1:\mathfrak g_1\longrightarrow \mathfrak
 g_1$ are linear maps
 and $d_2:\wedge^{3}\mathfrak g_0
\longrightarrow\mathfrak g_{1}$ is a completely skew-symmetric
trilinear map, and for any $x_i\in \mathfrak
g_0~(i=1,\cdot\cdot\cdot,5),~a\in \mathfrak
 g_1$, they satisfy the following conditions:
 $$d_0h=hd_1,\eqno(6.6)$$
\begin{eqnarray*}&&hd_2(x_1,x_2,x_3)+d_0l_3(x_1,x_2,x_3)\\&=&l_3(d_0x_1,x_2,x_3)+l_3(x_1,d_0x_2,x_3)+l_3(x_1,x_2,d_0x_3)
+\lambda(l_3(d_0x_1,d_0x_2,x_3)\\&&+l_3(x_1,d_0x_2,d_0x_3)+l_3(d_0x_1,x_2,d_0x_3))+\lambda^{2}l_3(d_0x_1,d_0x_2,d_0x_3)
,~~~~~~~~~~~~~~~~~~~~(6.7)
\\&&d_2(x_1,x_2,ha)+d_1l_3(x_1,x_2,a)\\&=&l_3(d_0x_1,x_2,a)+l_3(x_1,d_0x_2,a)+l_3(x_1,x_2,d_{1}a)
+\lambda(l_3(d_0x_1,d_0x_2,a)\\&&+l_3(x_1,d_0x_2,d_1a)+l_3(d_0x_1,x_2,d_1a))+\lambda^{2}l_3(d_0x_1,d_0x_2,d_1a)
,~~~~~~~~~~~~~~~~~~~~~~~~(6.8)
\\&&d_1l_5(x_1,x_2,x_3,x_4,x_5)-\sum_{k=1}^{5}\lambda^{k-1}l_{5}^{k}(x_1,x_2,x_3,x_4,x_5)
\\&=&d_2(x_3,l_3(x_1,x_2,x_4),x_5)+d_2(l_3(x_1,x_2,x_3),x_4,x_5)+d_2(x_3,x_4,l_3(x_1,x_2,x_5))
\\&&-d_2(x_1,x_2,l_3(x_3,x_4,x_5))
-l_3(x_1,x_2,d_2(x_3,x_4,x_5))
-\lambda(l_3(d_0x_1,x_2,d_2(x_3,x_4,x_5))\\&&+l_3(x_1,d_0x_2,d_2(x_3,x_4,x_5))+\lambda
l_3(d_0x_1,d_0x_2,d_2(x_3,x_4,x_5))) +l_3(x_3,d_2(x_1,x_2,x_4),x_5)
\\&&+\lambda(l_3(d_0x_3,d_2(x_1,x_2,x_4),x_5)+l_3(x_3,d_2(x_1,x_2,x_4),d_0x_5)\\&&+\lambda
l_3(d_0x_3,d_2(x_1,x_2,x_4),d_0x_5)) +l_3(d_2(x_1,x_2,x_3),x_4,x_5)
+\lambda(l_3(d_2(x_1,x_2,x_3),d_0x_4,x_5)\\&&+l_3(d_2(x_1,x_2,x_3),x_4,d_0x_5)+\lambda
l_3(d_2(x_1,x_2,x_3),d_0x_4,d_0x_5))
 +l_3(x_3,x_4,d_2(x_1,x_2,x_5))
\\&&+\lambda(l_3(d_0x_3,x_4,d_2(x_1,x_2,x_5))+l_3(x_3,d_0x_4,d_2(x_1,x_2,x_5))\\&&+\lambda
l_3(d_0x_3,d_0x_4,d_2(x_1,x_2,x_5))),~~~~~~~~~~~~~~~~~~~~~~~~~~~~~~~~~~~~~~~~~~~~~~~~~~~~~~~~~~~~(6.9)
 \end{eqnarray*}
where $l_{5}^{k} = l_{5}(\underbrace{I,\cdot\cdot\cdot ,
d_{0},\cdot\cdot\cdot , I}_{d_0~~ \hbox{appears} ~~k
~~\hbox{times}}),~ \hbox{where}~ I \hbox{~is ~identity~ map~ on}~
\mathfrak g_0.$ Denote it by $(\mathfrak g,d)$.\\

A differential 3-Lie 2-algebra is said to be skeletal (strict) if
$h=0$~($l_5=0,d_2=0$).\\

{\bf Proposition 6.2.} Let $(\mathfrak g,d)$ be a differential 3-Lie
2-algebra. We have

(i) If $(\mathfrak g,d)$ is skeletal or strict, then $(\mathfrak
g_0,[ \ , \ , \ ]_{\mathfrak g_0},l_3,d_0)$ is a differential 3-Lie
algebra, where $[ x ,y,z ]_{\mathfrak g_0}=l_3(x,y,z)$ for any
$x,y,z\in {\mathfrak g_0}.$

(ii) If $(\mathfrak g,d)$ is strict, then $(\mathfrak g_1,[ \ , \ ,
\ ]_{\mathfrak g_1},d_1)$ is a differential 3-Lie algebra, where
$$[a,b,c]_{\mathfrak
g_1}=l_3(h(a),h(b),c)=l_3(h(a),b,h(c))=l_3(a,h(b),h(c)),\eqno(6.10)$$
for any $a,b,c\in {\mathfrak g_1}.$

(iii) If $(\mathfrak g,d)$ is skeletal or strict, then $(\mathfrak
g_1,\rho,d_1)$ is a representation of $(\mathfrak g_0,[ \ , \ , \
]_{\mathfrak g_0},d_0)$ with bilinear map $\rho:\wedge^{2}{\mathfrak
g_0}\longrightarrow gl({\mathfrak g_1})$ given by
$\rho(x,y)a=l_3(x,y,a)$ for any $x,y\in {\mathfrak g_0},~a\in
{\mathfrak g_1}.$

\begin{proof} (i) can  be obtained by (6.3) and (6.7). (ii) can  be checked  by (6.3), (6.4) and (6.8).

 (iii) By (6.3) and (6.4), (2.3) and (2.4) hold.
Using (6.8),
\begin{eqnarray*}&&(\rho(d_0x,y)+\rho(x,d_0y)+\lambda\rho(d_0x,d_0y))(a+\lambda d_1a)
+\rho(x,y)d_1)a-d_1\rho(x,y)a
\\&=&l_3(d_0x,y,a)+l_3(x,d_0y,a)+\lambda l_3(d_0x,d_0y,a)+\lambda(l_3(d_0x,y,d_1a)+l_3(x,d_0y,d_1a)\\&&+\lambda l_3(d_0x,d_0y,d_1a))
+l_3(x,y,d_1a)-d_1l_3(x,y,a)
\\&=&0,\end{eqnarray*}
that is, (2.5) holds.
 Thus, $(\mathfrak g_1,\rho,d_1)$ is a representation of $(\mathfrak
g_0,[ \ , \ , \ ]_{\mathfrak g_0},d_0)$.

\end{proof}

{\bf Theorem 6.3.} There is a one-to-one correspondence between
skeletal differential 3-Lie 2-algebras and 3-cocycles of
differential 3-Lie algebras.
\begin{proof} Let $(\mathfrak g,d)$ be a skeletal differential 3-Lie
2-algebras. In view of Proposition 6.2,  we can consider the
cohomology of differential 3-Lie algebra $(\mathfrak g_0,[ \ , \ , \
]_{\mathfrak g_0},d_0)$ with coefficients in the representation
$(\mathfrak g_1,\rho,d_1)$ with $\rho$ defined in Proposition 6.2,
then $(l_5,d_2)$ is a 3-cocycle. In fact, for any
$X_i=x_{2i-1}\wedge x_{2i}\in \mathfrak g_0\wedge \mathfrak g_0
~(i=1,2,3),~ x_7\in \mathfrak g_0$,
$$\partial_{D_\lambda}(l_5,d_2)=(\partial l_5,\partial_\lambda
d_2+(-1)^{3}\delta l_5).$$
Combining (3.1) and (6.5), we have
\begin{eqnarray*}&&\partial l_5 (X_1,X_2,X_3,x_7)
\\&=&-l_5([X_1,X_2]_F,X_3,x_7)+l_5(X_1,[X_2,X_3]_F,x_7)-l_5(X_2,X_3,[X_1,x_7])+l_5(X_1,X_3,[X_2,x_7])\\&&-l_5(X_1,X_2,[X_3,x_7])
+\rho(X_1)l_5(X_2,X_3,x_7)-\rho(X_2)l_5(X_1,X_3,x_7)+\rho(X_3)l_5(X_1,X_2,x_7)\\&&+\rho(x_6,x_7)l_5(X_1,X_2,x_5)+\rho(x_7,x_5)l_5(X_1,X_2,x_6)
\\&=&0.\end{eqnarray*}
By (2.7), (3.2) and (6.9), there holds that
\begin{eqnarray*}&&(\partial_\lambda d_2-\delta l_5)
(x_1,x_2,x_3,x_4,x_5)
\\&=&-d_2(x_3,[x_1,x_2,x_4]_{\mathfrak g_0},x_5)-d_2([x_1,x_2,x_3]_{\mathfrak g_0},x_4,x_5)-d_2(x_3,x_4,[x_1,x_2,x_5]_{\mathfrak g_0})
\\&&+d_2(x_1,x_2,[x_3,x_4,x_5]_{\mathfrak g_0})+\hat{\rho}(x_1,x_2)d_2(x_3,x_4,x_5)-\hat{\rho}(x_3,x_4)d_2(x_1,x_2,x_5)
\\&&-\hat{\rho}(x_4,x_5)d_2(x_1,x_2,x_3)-\hat{\rho}(x_5,x_3)d_2(x_1,x_2,x_4)
\\&&-\sum_{k=1}^{5}\lambda^{k-1}l_{5}^{k}(x_1,x_2,x_3,x_4,x_5)+d_1l_5
(x_1,x_2,x_3,x_4,x_5)\\&=&0.
 \end{eqnarray*}
 Hence,
$(l_5,d_2)$ is a 3-cocycle.

Conversely, suppose that $(\varphi,\psi)$ is a 3-cocycle of
$(\mathfrak g_0,l_3,d_0)$ with coefficients in the representation
$(V,\rho,d_V)$. Then $(\mathfrak g,d)$ is a skeletal differential
3-Lie 2-algebra, where $d=(d_0,d_1=d_V,d_2=\psi)$
 and $\mathfrak g=(\mathfrak g_0,\mathfrak g_1=V,h=0,l_3,l_5=\varphi)$
 with $l_3(x,y,v)=\rho(x,y)v$ for any $x,y\in \mathfrak g_0,~v\in
 V$. We omit the details.
\end{proof}

{\bf Definition 6.4. } A crossed module of differential 3-Lie
2-algebra $(\mathfrak g,d)$ is a sextuple $((\mathfrak g_0, [ \ , \
, \ ]_{\mathfrak g_0}),(\mathfrak g_1, [ \ , \ , \ ]_{\mathfrak
g_1}), h,\rho,d_0,d_1)$, where $(\mathfrak g_0, [ \ , \ , \
]_{\mathfrak g_0},d_0)$ and $(\mathfrak g_1, [ \ , \ , \
]_{\mathfrak g_1},d_1)$ are differential 3-Lie algebras,
$h:\mathfrak g_1\longrightarrow \mathfrak g_0$ is a homomorphism of
 differential 3-Lie algebras, and $(\mathfrak g_1,\rho,d_1)$ is a representation
of $(\mathfrak g_0, [ \ , \ , \ ]_{\mathfrak g_0},d_0)$, and the
followings are satisfied:
$$h\rho(x,y)a=[x,y,h(a)]_{\mathfrak g_0},\eqno(6.11)$$
$$\rho(h(a),h(b))c=[a,b,c]_{\mathfrak g_1},~~\rho(x,h(a))b=-\rho(x,h(b))a,\eqno(6.12)$$
for any $x,y\in {\mathfrak g_0},~a,b,c\in {\mathfrak g_1}.$\\

{\bf Theorem 6.5.} There is a one-to-one correspondence between
strict differential 3-Lie 2-algebras and crossed modules of
differential 3-Lie 2-algebras.

\begin{proof} On the one hand, let $(\mathfrak g,d)$ be a strict differential 3-Lie 2-algebra of
weight $\lambda$, we show that $((\mathfrak g_0, [ \ , \ , \
]_{\mathfrak g_0}),(\mathfrak g_1, [ \ , \ , \ ]_{\mathfrak g_1}),
h,\rho,d_0,d_1)$ is a crossed module of $(\mathfrak g,d)$ with
notions as in Proposition 6.2. Based on Proposition 6.2, we only
need to check that (6.11), (6.12) hold and $h$ is a homomorphism of
differential 3-Lie algebras.

Clearly, (6.11) holds by (6.1). In view of (6.10) and (6.1), (6.12)
holds. For any $a,b,c\in {\mathfrak g_1}$,
$$h([ a, b, c ]_{\mathfrak
g_1})=hl_3(h(a),h(b),c)=l_3(h(a),h(b),h(c))=[
h(a),h(b),h(c)]_{\mathfrak g_0},$$ and $d_0h=hd_1$ by (6.6). We
induce that $h$ is a homomorphism of differential 3-Lie algebras.

On the other hand, it can be verified analogously.

\end{proof}

\section{Appendix: Proof of Proposition 3.3}
\begin{proof} Define $f_{i,j}^{k}$ by
$$f_{i,j}^{k} = f(\underbrace{I,\cdot\cdot\cdot , d_{\mathfrak
g},\cdot\cdot\cdot , I}_{d_{\mathfrak g}~~ \hbox{appears} ~~k
~~\hbox{times  but not in the ith,jth places}}),$$ and $f_{j}^{k}$
is just the $f_{i,j}^{k} $ for no $i$, where $I$ is the identity map
on $\mathfrak g$.

 For any $f\in C_{3-{\hbox{Lie}}}^{n}(\mathfrak
g,V)$ and $X_1,\cdot\cdot\cdot,X_{n}\in \mathfrak g\wedge \mathfrak
g,z\in \mathfrak g$, suppose $X_i=x_i\wedge y_i$, and denote
\begin{eqnarray*}\tilde{X_{11}}=[x_j,y_j,
x_i]\wedge y_i+x_i\wedge [x_j,y_j,y_i],
\end{eqnarray*}
\begin{eqnarray*}\tilde{X_{12}}&=&([d_{\mathfrak g}(x_j),y_j,x_i]\wedge
y_i+x_i\wedge [x_j,y_j,y_i])+([x_j,d_{\mathfrak g}(y_j),x_i]\wedge
y_i+x_i\wedge [x_j,d_{\mathfrak g}(y_j),y_i])\\&&+
([x_j,y_j,d_{\mathfrak g}(x_i)]\wedge y_i+d_{\mathfrak g}(x_i)\wedge
[x_j,y_j,y_i]),
\end{eqnarray*}
\begin{eqnarray*}\tilde{X_{13}}&=& ([d_{\mathfrak g}(x_j),d_{\mathfrak
g}(y_j),x_i]\wedge y_i+x_i\wedge [d_{\mathfrak g}(x_j),d_{\mathfrak
g}(y_j),y_i]+([d_{\mathfrak g}(x_j),y_j,d_{\mathfrak g}(x_i)]\wedge
y_i\\&&+d_{\mathfrak g}(x_i)\wedge [d_{\mathfrak g}(x_j), y_j,y_i])+
([d_{\mathfrak g}(x_j),y_j,x_i]\wedge d_{\mathfrak g}(y_i)+x_i\wedge
[d_{\mathfrak g}(x_j),y_j,d_{\mathfrak g}(y_i)])\\&&+
([x_j,d_{\mathfrak g}(y_j),d_{\mathfrak g}(x_i)]\wedge
y_i+d_{\mathfrak g}(x_i)\wedge [ x_j,d_{\mathfrak g}(y_j),y_i]) +
([x_j,d_{\mathfrak g}(y_j),x_i]\wedge d_{\mathfrak
g}(y_i)\\&&+x_i\wedge [ x_j,d_{\mathfrak g}(y_j),d_{\mathfrak
g}(y_i)]) + ([x_j,y_j,d_{\mathfrak g}(x_i)]\wedge d_{\mathfrak
g}(y_i)+d_{\mathfrak g}(x_i)\wedge [ x_j,y_j,d_{\mathfrak g}(y_i)]),
\end{eqnarray*}
\begin{eqnarray*}\tilde{X_{14}}&=&([d_{\mathfrak g}(x_j),d_{\mathfrak g}(y_j),d_{\mathfrak
g}(x_i)]\wedge y_i+d_{\mathfrak g}(x_i)\wedge [d_{\mathfrak
g}(x_j),d_{\mathfrak g}(y_j),y_i] \\&&+([d_{\mathfrak
g}(x_j),y_j,x_i]\wedge d_{\mathfrak g}(y_i)+x_i\wedge [d_{\mathfrak
g}(x_j), d_{\mathfrak g}(y_j),d_{\mathfrak g}(y_i)]) \\&&+
([x_j,d_{\mathfrak g}(y_j),d_{\mathfrak g}(x_i)]\wedge d_{\mathfrak
g}(y_i)+d_{\mathfrak g}(x_i)\wedge [x_j,d_{\mathfrak
g}(y_j),d_{\mathfrak g}(y_i)]),
\end{eqnarray*}
\begin{eqnarray*}\tilde{X_{15}}=[d_{\mathfrak g}(x_j),d_{\mathfrak g}(y_j),d_{\mathfrak
g}(x_i)]\wedge d_{\mathfrak g}(y_i)+d_{\mathfrak g}(x_i)\wedge
[d_{\mathfrak g}(x_j),d_{\mathfrak g}(y_j),d_{\mathfrak g}(y_i)],
\end{eqnarray*}
$$\tilde{X_{16}}=[x_j,y_j,z],~\tilde{X_{17}}=[d_{\mathfrak
g}(x_j),y_j,z]+ [x_j,d_{\mathfrak g}(y_j),z]+ [x_j,y_j,d_{\mathfrak
g}(z)],$$ and
$$\tilde{X_{18}}=[d_{\mathfrak g}(x_j),d_{\mathfrak g}(y_j),z]+
[d_{\mathfrak g}(x_j),y_j,d_{\mathfrak g}(z)]+ [x_j,d_{\mathfrak
g}(y_j),d_{\mathfrak g}(z)],~\tilde{X_{19}}= [d_{\mathfrak
g}(x_j),d_{\mathfrak g}(y_j),d_{\mathfrak g}(z)].
$$
On the one hand, we calculate $\delta\partial(f)$.
\begin{eqnarray*}&&\delta\partial(f)(X_1,\cdot\cdot\cdot,X_{n},z)\\&=&\sum_{k=1}^{2n+1}\lambda^{k-1}\partial
f^{k}(X_1,\cdot\cdot\cdot,X_{n},z)-d_{V}\partial(f)(X_1,\cdot\cdot\cdot,X_{n},z)\\&=&A+B+C+D+E-K,
\end{eqnarray*}
where
\begin{eqnarray*}A&=&\sum_{k=1}^{2n-3}\lambda^{k-1}\sum_{1\leq
j<i\leq n}(-1)^{j}
f_{2j-1,2j}^{k}(X_1,\cdot\cdot\cdot,X_{j-1},\tilde{X_{11}},X_{j+1},\cdot\cdot\cdot,X_n,z)
\\&&+\sum_{k=1}^{2n-2}\lambda^{k-1}\sum_{1\leq j<i\leq n}(-1)^{j}
f_{2j-1,2j}^{k-1}(X_1,\cdot\cdot\cdot,X_{j-1},\tilde{X_{12}},X_{j+1},\cdot\cdot\cdot,X_n,z)
\\&&+\sum_{k=2}^{2n-1}\lambda^{k-1}\sum_{1\leq j<i\leq n}(-1)^{j}
f_{2j-1,2j}^{k-2}(X_1,\cdot\cdot\cdot,X_{j-1},\tilde{X_{13}},X_{j+1},\cdot\cdot\cdot,X_n,z)
\\&&+\sum_{k=3}^{2n}\lambda^{k-1}\sum_{1\leq j<i\leq n}(-1)^{j}
f_{2j-1,2j}^{k-3}(X_1,\cdot\cdot\cdot,X_{j-1},\tilde{X_{14}},X_{j+1},\cdot\cdot\cdot,X_n,z)
\\&&+\sum_{k=4}^{2n+1}\lambda^{k-1}\sum_{1\leq j<i\leq n}(-1)^{j}
f_{2j-1,2j}^{k-4}(X_1,\cdot\cdot\cdot,X_{j-1},\tilde{X_{15}},X_{j+1},\cdot\cdot\cdot,X_n,z),
\end{eqnarray*}
\begin{eqnarray*}B&=&\sum_{k=1}^{2n-2}\lambda^{k-1}\sum_{j=1}^{n}(-1)^{j}f_{2n-1}^{k}
(X_1,\cdot\cdot\cdot,X_{j-1},X_{j+1},\cdot\cdot\cdot,X_{n},\tilde{X_{16}})
\\&&+\sum_{k=1}^{2n-1}\lambda^{k-1}\sum_{j=1}^{n}(-1)^{j}f_{2n-1}^{k-1}(X_1,\cdot\cdot\cdot,X_{j-1},X_{j+1},\cdot\cdot\cdot,X_{n},\tilde{X_{17}})
\\&&+\sum_{k=2}^{2n}\lambda^{k-1}\sum_{j=1}^{n}(-1)^{j}f_{2n-1}^{k-2}(X_1,\cdot\cdot\cdot,X_{j-1},X_{j+1},\cdot\cdot\cdot,X_{n},\tilde{X_{18}})
\\&&+\sum_{k=3}^{2n+1}\lambda^{k-1}\sum_{j=1}^{n}(-1)^{j}f_{2n-1}^{k-3}(X_1,\cdot\cdot\cdot,X_{j-1},X_{j+1},\cdot\cdot\cdot,X_{n},\tilde{X_{19}}),
\end{eqnarray*}
\begin{eqnarray*}C&=&\sum_{k=1}^{2n-1}\lambda^{k-1}\sum_{j=1}^{n}(-1)^{j+1}\rho(x_j,y_j)
f^{k}(X_1,\cdot\cdot\cdot,X_{j-1},X_{j+1},\cdot\cdot\cdot,X_{n},z)\\&&+
\sum_{k=1}^{2n}\lambda^{k-1}\sum_{j=1}^{n}(-1)^{j+1}(\rho(d_{\mathfrak
g}(x_j),y_j)+\rho(x_j,d_{\mathfrak g}(y_j))
f^{k-1}(X_1,\cdot\cdot\cdot,X_{j-1},X_{j+1},\cdot\cdot\cdot,X_{n},z)
\\&&+\sum_{k=2}^{2n+1}\lambda^{k-1}\sum_{j=1}^{n}(-1)^{j+1}\rho(d_{\mathfrak
g}(x_j),d_{\mathfrak g}(y_j))
f^{k-2}(X_1,\cdot\cdot\cdot,X_{j-1},X_{j+1},\cdot\cdot\cdot,X_{n},z)
\\&=&
\sum_{k=1}^{2n-1}\lambda^{k-1}\sum_{j=1}^{n}(-1)^{j+1}\rho(x_j,y_j)
f^{k}(X_1,\cdot\cdot\cdot,X_{j-1},X_{j+1},\cdot\cdot\cdot,X_{n},z)\\&&+
\sum_{k=0}^{2n-1}\lambda^{k}\sum_{j=1}^{n}(-1)^{j+1}(\rho(d_{\mathfrak
g}(x_j),y_j)+\rho(x_j,d_{\mathfrak g}(y_j))
f^{k}(X_1,\cdot\cdot\cdot,X_{j-1},X_{j+1},\cdot\cdot\cdot,X_{n},z)
\\&&+\sum_{k=0}^{2n-1}\lambda^{k+1}\sum_{j=1}^{n}(-1)^{j+1}\rho(d_{\mathfrak
g}(x_j),d_{\mathfrak g}(y_j))
f^{k}(X_1,\cdot\cdot\cdot,X_{j-1},X_{j+1},\cdot\cdot\cdot,X_{n},z),
\end{eqnarray*}
\begin{eqnarray*}D&=&\sum_{k=1}^{2n-1}\lambda^{k-1}(-1)^{n+1}\rho(y_n,z)f^{k}(X_1,\cdot\cdot\cdot,X_{n-1},x_{n})
\\&&+\sum_{k=1}^{2n}\lambda^{k-1}(-1)^{n+1}\rho(d_{\mathfrak
g}(y_n),z)f^{k-1}(X_1,\cdot\cdot\cdot,X_{n-1},x_{n})
\\&&+\sum_{k=1}^{2n}\lambda^{k-1}(-1)^{n+1}\rho(y_n,d_{\mathfrak
g}(z))f^{k-1}(X_1,\cdot\cdot\cdot,X_{n-1},x_{n})
\\&&+\sum_{k=2}^{2n+1}\lambda^{k-1}(-1)^{n+1}\rho(d_{\mathfrak
g}(y_n),d_{\mathfrak
g}(z))f^{k-2}(X_1,\cdot\cdot\cdot,X_{n-1},x_{n}),
\end{eqnarray*}
\begin{eqnarray*}E&=&\sum_{k=1}^{2n-1}\lambda^{k-1}(-1)^{n+1}\rho(z,x_n)f^{k}(X_1,\cdot\cdot\cdot,X_{n-1},y_{n})
\\&&+\sum_{k=1}^{2n}\lambda^{k-1}(-1)^{n+1}\rho(d_{\mathfrak
g}(z),x_n)f^{k-1}(X_1,\cdot\cdot\cdot,X_{n-1},y_{n})
\\&&+\sum_{k=1}^{2n}\lambda^{k-1}(-1)^{n+1}\rho(z,d_{\mathfrak
g}(x_n))f^{k-1}(X_1,\cdot\cdot\cdot,X_{n-1},y_{n})
\\&&+\sum_{k=2}^{2n+1}\lambda^{k-1}(-1)^{n+1}\rho(d_{\mathfrak
g}(z),d_{\mathfrak
g}(x_n))f^{k-2}(X_1,\cdot\cdot\cdot,X_{n-1},y_{n}),
\end{eqnarray*}
and $$K=d_V(\partial f)(X_1,\cdot\cdot\cdot,X_{n},z)=K_1+K_2+K_3+K_4
$$
with
$$K_1= d_V\sum_{1\leq j<i\leq
n}(-1)^{j}f(X_1,\cdot\cdot\cdot,X_{j-1},X_{j+1},\cdot\cdot\cdot,X_{i-1},[X_j,X_i]_F,X_{i+1},\cdot\cdot\cdot,X_{n}
,z),$$
$$K_2=d_V\sum_{j=1}^{n}(-1)^{j}f(X_1,\cdot\cdot\cdot,\hat{X_j},\cdot\cdot\cdot,X_{n},[X_j,z]),$$
\begin{eqnarray*}K_3&=&d_V\sum_{j=1}^{n}(-1)^{j+1}\rho(X_j)f(X_1,\cdot\cdot\cdot,X_{j-1},X_{j+1},\cdot\cdot\cdot,X_{n},z)\\
&=&\sum_{j=1}^{n}(-1)^{j+1} (\rho(x_j,y_j)d_{V}+\rho(d_{\mathfrak
g}(x_j),y_j)+\rho(x_j,d_{\mathfrak g}(y_j))+\lambda\rho(d_{\mathfrak
g}(x_j),d_{\mathfrak g}(y_j))d_{V}\\&&+\lambda\rho(d_{\mathfrak
g}(x_j),y_j)d_{V}+\lambda\rho(x_j,d_{\mathfrak
g}(y_j))d_{V}+\lambda^2\rho(d_{\mathfrak g}(x_j),d_{\mathfrak
g}(y_j))d_{V})\\&&\times
f(X_1,\cdot\cdot\cdot,X_{j-1},X_{j+1},\cdot\cdot\cdot,X_{n},z)
,\end{eqnarray*} and
$$K_4=(-1)^{n+1}d_{V}(\rho(y_{n},z)f(X_1,\cdot\cdot\cdot,X_{n-1},x_{n})+\rho(z,x_{n})f(X_1,\cdot\cdot\cdot,X_{n-1},y_{n}))).$$
On the other hand,
$$\partial_{\lambda}\delta(f)(X_1,\cdot\cdot\cdot,X_{n},z)
=F+G+H+J,$$ where
\begin{eqnarray*}F&=&
\partial_{\lambda}\delta(f)(X_1,\cdot\cdot\cdot,X_{n},z)
\\&=&\sum_{1\leq j<
i\leq {n}} (-1)^{j}\delta
f(X_1,\cdot\cdot\cdot,X_{j-1},X_{j+1},\cdot\cdot\cdot,X_{i-1},[X_j,X_i]_F,X_{i+1},\cdot\cdot\cdot,X_{n}
,z)\\&=& \sum_{1\leq j< i\leq {n}}
(-1)^{j}\sum_{k=1}^{2n-1}\lambda^{k-1}f^{k}(X_1,\cdot\cdot\cdot,X_{j-1},X_{j+1},\cdot\cdot\cdot,X_{i-1},[X_j,X_i]_F,X_{i+1},\cdot\cdot\cdot,X_{n}
,z)\\&&- \sum_{1\leq j< i\leq {n}}
(-1)^{j}d_Vf(X_1,\cdot\cdot\cdot,X_{j-1},X_{j+1},\cdot\cdot\cdot,X_{i-1},[X_j,X_i]_F,X_{i+1},\cdot\cdot\cdot,X_{n}
,z)\\&=&F_1-F_2,
\end{eqnarray*}
\begin{eqnarray*}G&=&\sum_{j=1}^{n}(-1)^{j}\delta
f(X_1,\cdot\cdot\cdot,X_{j-1},X_{j+1},\cdot\cdot\cdot,X_{n},[X_j,z])
\\&=&\sum_{j=1}^{n}(-1)^{j}\sum_{k=1}^{2n-1}\lambda^{k-1}f^{k}(X_1,\cdot\cdot\cdot,X_{j-1},X_{j+1},\cdot\cdot\cdot,X_{n},[X_j,z])
\\&&-\sum_{j=1}^{n}(-1)^{j}d_Vf(X_1,\cdot\cdot\cdot,X_{j-1},X_{j+1},\cdot\cdot\cdot,X_{n},[X_j,z])\\&=&G_1-G_2,
\end{eqnarray*}
According to (2.7),
\begin{eqnarray*}H&=&\sum_{j=1}^{n}(-1)^{j+1}\hat{\rho}(X_j)\delta
f(X_1,\cdot\cdot\cdot,X_{j-1},X_{j+1},\cdot\cdot\cdot,X_{n},z)\\&=&\sum_{j=1}^{n}(-1)^{j+1}
(\rho(X_j)+\lambda\rho(d_{\mathfrak
g}(x_j),y_j)+\lambda\rho(x_j,d_{\mathfrak
g}(y_j))+\lambda^2\rho(d_{\mathfrak g}(x_j),d_{\mathfrak
g}(y_j)))\\&&\times
(\sum_{k=1}^{2n-1}\lambda^{k-1}f^{k}(X_1,\cdot\cdot\cdot,X_{j-1},X_{j+1},\cdot\cdot\cdot,X_{n},z)-d_{V}f(X_1,\cdot\cdot\cdot,X_{j-1},X_{j+1},\cdot\cdot\cdot,X_{n},z)
\\&=&H_1-H_2,\end{eqnarray*}
and
\begin{eqnarray*}J&=&(-1)^{n+1}(\hat{\rho}(y_{n},z)\delta
f(X_1,\cdot\cdot\cdot,X_{n-1},x_{n})+\hat{\rho}(z,x_{n})\delta
f(X_1,\cdot\cdot\cdot,X_{n-1},y_{n})).
\end{eqnarray*}
Clearly, we have $$A-K_1=F_1-F_2,B-K_2=G_1-G_2,C-K_3=H_1-H_2.$$
Similarly to the case of $C-K_3=H_1-H_2$, we can prove that
$D+E-K_4=J$. Thus, $\partial_\lambda\delta=\delta\partial$.
\end{proof}


\begin{center}{\textbf{Acknowledgments}}
\end{center}
Project supported by the National Natural Science Foundation of
China (No. 11871421, No. 11401530) and the Natural Science
Foundation of Zhejiang Province of China ( No. LY19A010001 ).

\begin{center} {\textbf{Data availability}}
\end{center}
All datasets underlying the conclusions of the paper are available
to readers.


\end{document}